\documentclass[11pt]{article}

\usepackage[utf8]{inputenc}

\usepackage{graphicx}
\usepackage{url}
\usepackage{amsmath}
\usepackage{amsfonts}
\usepackage{verbatim}
\usepackage{amssymb}
\usepackage{amsthm}
\usepackage{color}
\usepackage{float}
\usepackage[english]{babel}
\usepackage[autostyle]{csquotes}
\usepackage{mathtools}
\usepackage[usenames,dvipsnames]{xcolor}
\usepackage{mdwlist}
\usepackage{caption}
\usepackage{subcaption}
\usepackage{setspace}
\usepackage{bbm}
\usepackage{wrapfig}
\usepackage{empheq}
\usepackage{hyperref}
\hypersetup{
    colorlinks=true,
    linkcolor=blue,
    filecolor=magenta,      
    urlcolor=blue,
    citecolor=blue,
}

\usepackage[top=1in, bottom=1in, left=1in, right=1in]{geometry}

\setcaptionmargin{0.25in}

\title{Reflections in excitable media linked to existence and stability of one-dimensional spiral waves}
\author{Stephanie Dodson\thanks{Department of Mathematics, University of California, Davis, Davis, CA: \texttt{sadodson@ucdavis.edu}}
\and Timothy J. Lewis\thanks{Department of Mathematics, University of California, Davis, Davis, CA: \texttt{tjlewis@ucdavis.edu}}}
\date{\today}

\begin{document}

\maketitle

% REQUIRED
\begin{abstract}
 \noindent When propagated action potentials in cardiac tissue interact with local heterogeneities,  reflected waves can sometimes be induced. These reflected waves have been associated with the onset of cardiac arrhythmias, and while their generation is not well understood, their existence is linked to that of one-dimensional (1D) spiral waves. Thus, understanding the existence and stability of 1D spirals plays a crucial role in determining the likelihood of the unwanted reflected pulses. Mathematically, we probe these issues by viewing the 1D spiral as a time-periodic antisymmetric source defect. Through a combination of direct numerical simulation and continuation methods, we investigate existence and stability of a 1D spiral wave in a qualitative ionic model to determine how the systems propensity for reflections are influenced by system parameters. Our results support and extend a previous hypothesis that the 1D spiral is an unstable periodic orbit that emerges through a global rearrangement of heteroclinic orbits and we identify key parameters and physiological processes that promote and deter reflection behavior.  
\end{abstract}

% REQUIRED
\textbf{Keywords:} one-dimensional spiral waves, source-defects, wave reflections, cardiac arrhythmia, Morris-Lecar

% REQUIRED
\textbf{AMS subject classifications:} 35B36, 35K57, 35P30, 92-10

%35B36: Pattern Formation
% 35K57: Reaction-diffusion equations
% 92-10: Mathematical modeling or simulation for problems pertaining to biology
% 35P30: Nonlinear eigenvalue problems and nonlinear spectral theory for PDEs

\section{Introduction}

The ability of the heart to effectively pump blood to the body relies on the coordinated spread of electrochemical waves (propagated action potentials) through the cardiac tissue. If the spread of the waves are disrupted, life-threatening arrhythmias may arise. The most severe form of arrhythmias involve reentry behavior, which is composed of self-perpetuating circulating waves that underlies arrhythmias such as ventricular tachycardia and ventricular fibrillation. Reentrant arrhythmias are thought to be induced by propagated action potentials interacting with functional or structural heterogeneities in the cardiac tissue that lead to symmetry breaking of the wave via local one-way propagation block.

One of the mechanisms that is thought to induce the onset of reentrant arrhythmias is reflection \cite{aj80,vsjjg88,ab11,kb15}. Usually, when a propagated cardiac action potential interacts with a local tissue heterogeneity, the action potential is either blocked and annihilated or it successfully propagates across the heterogeneity. However, in some cases, action potentials can successfully cross the heterogeneity and then trigger a reflected action potential, i.e., an action potential that propagates only in the retrograde direction. A locally reflected pulse breaks the symmetry of the wave of activation spreading through the heart and thus can initiate reentry. Reflection has been directly observed in cardiac Purkinje fiber \cite{wch72, aj80, jm81, a83}, ventricular tissue \cite{am81, rjm84, da85}, and atrial tissue \cite{la89} in experimental preparations with regions of depressed excitability or reduced electrical coupling. More recently, reflection was observed in monolayers of ventricular cardiomyocytes with sites of abrupt tissue expansions \cite{Auerbach:2011jm}. Action potential reflection has also been demonstrated in computational models of cardiac tissue \cite{cb92, l98, Auerbach:2011jm, kb15}. Furthermore, reflections have been observed in experimental preparations of neural tissue at the axon level \cite{hcl76, Baccus:1998hh, bbsm00} and the network level \cite{cc89, xhtw07}, as well as corresponding neural models \cite{gr74, zb94, ge11}.

 Previous work began to uncover the dynamical mechanisms underlying reflections of pulses. In 1996, Ermentrout \& Rinzel conducted numerical experiments on a cable model of an axon using Morris-Lecar dynamics and simulated the effects of an abrupt change in the diameter on action potential propagation \cite{er96}. When the diameter of the distal portion of the axon was sufficiently large, the action potential failed to propagate across the heterogeneity; whereas for sufficiently small distal axon diameters, the action potential successfully propagated across the heterogeneity. At intermediate diameters, the action potential propagated across the heterogeneity and then gave rise to a reflected pulse. Moreover, as the diameter was tuned to be closer and closer to a critical value, more and more reflections seemed to occur with successive counter-propagating action potentials. These simulations suggested that an unstable periodic orbit underlies the reflection behavior, and that the closer that the incoming wave brings the system to the periodic orbit, the more reflections occur. The authors pointed out that the unstable periodic orbit was a so-called one-dimensional (1D) spiral wave that was first described by Kopell \& Howard in \cite{kh81}. 1D spiral waves are `source defects’ \cite{ss04} that consist of a non-excited core that sheds anti-phase counter-propagating pulses.

\begin{figure}[ht]
    \centering
        \includegraphics[width=0.8\textwidth]{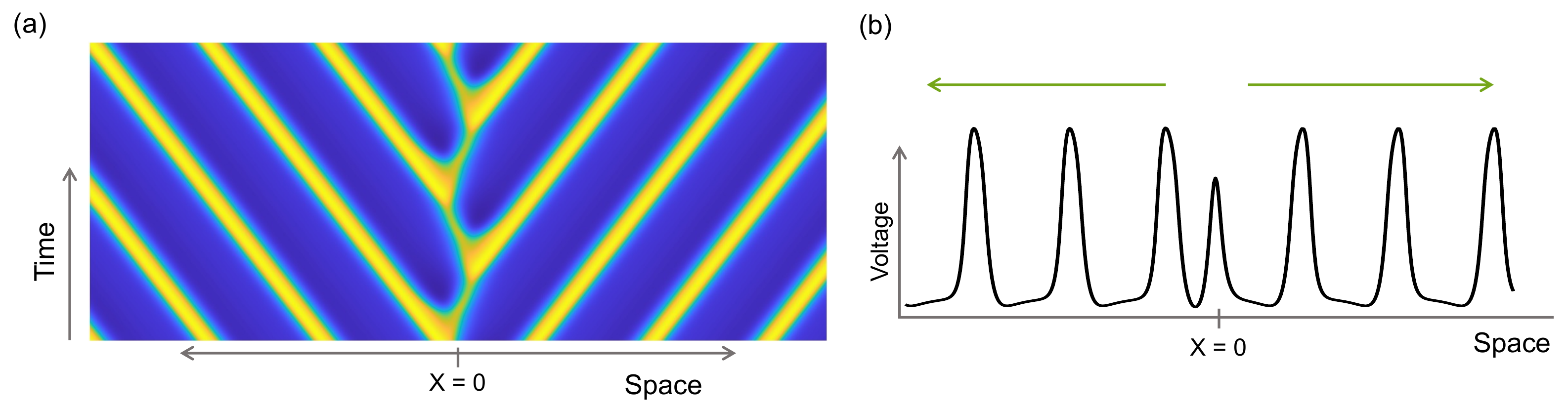}
        \caption{ (a) Space-time plot of a 1D spiral wave. (b) Example of spatial profile in 1D spiral. Green arrows above peaks provide direction of motion in the far-field. } \label{fig:reflection_example}
\end{figure}

An example of a series of reflections is given in Figure~\ref{fig:reflection_example}. The pattern has been dubbed a `1D spiral wave’ because its spatial structure resembles a radial slice through a two-dimensional rotating spiral wave. Waves emerge from the core on alternating sides and arrange into regular periodic structures as they propagate into the far-field. However, we emphasize that, despite the similarities in the names, the 1D and 2D spirals are different mathematical structures. Our focus here is specifically on the 1D spiral wave. 

To probe the link between reflection and the 1D spiral, Cytrynbaum \& Lewis \cite{cl09}  considered a homogeneous Morris-Lecar cable model and examined the interaction of a propagated action potential with a transient refractory region, i.e., a functional heterogeneity. By tuning the height of the refractory bump, they observed similar behavior as Ermentrout and Rinzel in \cite{er96}. However, reflection behavior was only seen when the excitability parameter $\epsilon$ was below a critical value $\epsilon_*$. Above this value $(\epsilon > \epsilon_*$), no evidence for the existence of the 1D spiral was detected, and no signs of reflections were observed (see Figure~\ref{fig:bifur_diag}). Numerical experiments near $\epsilon_*$ lead them to postulate that the unstable periodic orbit is created through a heteroclinic bifurcation at $\epsilon_*$. 

Pulse generators that resemble a one-sided 1D spiral have been analyzed in \cite{Prat:2005bp,Nishiura:2007fe,Teramoto:2009hd,Yadome:2014gt}. These structures were created by a jump-type heterogeneity and were also found to originate via a heteroclinic bifurcation \cite{Yadome:2014gt}. However, there are distinct differences from the scenario of the 1D spiral presented here. Specifically, the pulse generators of \cite{Yadome:2014gt} are stable and produce a variety of heterogeneity-induced patterns that have not been linked to reflections. 

Understanding how system parameters influence the existence of 1D spirals and the ability of the system to generate reflections is vitally important for identifying the physiological mechanisms that contribute to reflection-mediated arrhythmias. Even though the 1D spiral structure is unstable, the existence of the periodic orbit allows for new spatiotemporal dynamics, in that initial conditions near it lead to a transient series of spontaneous reflections. The specific physiological conditions that promote or deter reflections remain poorly understood, and the dynamical mechanisms creating the periodic orbit are unverified. It has been suggested \cite{er96} that reflection phenomena occur more robustly when local dynamics exhibit type I excitability (excitability that is associated with saddle-type threshold behavior) than type II excitability (quasi-threshold behavior) \cite{re89}. 

Here, we further probe these questions in the Morris-Lecar model. To do so, we interpret the system's propensity for reflections in two ways: (1) how conditions alter the probability of reflections occurring when the 1D spiral exists, and (2) how conditions modify the 1D spiral existence window. Previous studies detected the 1D spiral with direct numerical simulations, but these methods are inefficient for finding unstable solutions. Therefore, to investigate the 1D spiral existence and properties, we develop more refined numerical techniques that solve directly for the 1D spiral as a spatiotemporal equilibrium solution. Numerical methods used and developed will be outlined in the main text, and full implementation details are included in the Appendix. 

This paper is organized as follows. Section~\ref{sec:prelims} describes the Morris-Lecar model and introduces the mathematical framework of traveling pulses, periodic waves, reflections, and appearance of 1D spirals. In Section~\ref{sec:eigenvalues_likelihood}, we investigate scenarios when the 1D spiral does exist and understand how model parameters influence the likelihood of reflections by linking the stability of the 1D spiral with the probability of generating reflections. Conditions that alter the existence of 1D spirals and reflections are probed in Section~\ref{sec:1dspiral_existence}. We further verify the proposed heteroclinic structure of the 1D spiral bifurcation by confirming there are no additional local bifurcations of the slow pulse and accurately computing the period of the orbit near the bifurcation point. To compute the 1D spiral existence interval and it’s dependence on parameters, we define novel computational techniques to directly compute the heteroclinic bifurcation point. Our results are interpreted in the context of the local excitability type in Section~\ref{sec:typeI_typeII} to test the hypothesis that type I local dynamics promote reflections. Finally, in Section~\ref{sec:discussion} the stability and existence results are combined with the physiological mechanisms to discuss which processes could be most involved in generating reflections and provide predictions on how to decrease reflections in more complex systems. 

%%%%%% %%%%%% %%%%%% %%%%%% %%%%%% %%%%%% %%%%%% %%%%%% %%%%%% 
%% MATHEMATICAL FORMULATION %%%%%%%
%%%%%% %%%%%% %%%%%% %%%%%% %%%%%% %%%%%% %%%%%% %%%%%% %%%%%% 
\section{Mathematical Formulation} \label{sec:prelims}

Our primary goal is to understand conditions that lead to the existence of the unstable 1D spiral and impact its stability in the Morris-Lecar model. Previous studies relied on numerical simulations and shooting methods that are impractical and inefficient at finding unstable solutions. Our computational methods expand on techniques used to compute traveling waves and patterns as equilibrium solutions. Therefore, after a brief description of the Morris-Lecar system and parameters, we provide a short description of traveling waves in reaction-diffusion equations and the computational methods for these systems. We then end the section with a preliminary description of the appearance of the 1D spiral.

\subsection{Reaction-diffusion model with Morris-Lecar dynamics}

Reaction-diffusion equations in one spatial dimension take the form 
\begin{align} \label{eqn:rxn_diff}
U_t = D U_{xx} + F(U; \mu), \ \ U \in \mathbb{R}^n, \ D \in \mathbb{R}^{n \times n}, \ \mu \in \mathbb{R}^p, \ x \in \mathbb{R}, \ t \in \mathbb{R}
\end{align}
where $U = (u_1, u_2, \cdots, u_n)^T$ is a vector of species that diffuse at non-negative rates given by elements of the diagonal matrix $D$. The diffusion term couples nearby areas in space and the nonlinear reaction term $F(U;\mu)$ defines local excitable kinematic reactions, where $\mu$ represents the $p$ system parameters.

The Morris-Lecar model \cite{ml81} is a two-variable canonical model of excitable media, with variables $U = (V,n)^T$ representing the membrane potential $V$ (in mV) and the potassium channel gating variable $n$. Diffusive coupling is only throught the membrane potential. While idealized, the model is still framed in terms of biophysical variables and ionic currents, and parameters can be tuned so that the it exhibits a wide range of phenomena observed in cardiac and neural tissue. Thus, it is flexible and relatively easy to analyze, and results can be tied to biophysical mechanisms. 

Because we will be investigating the parameter impacts on the likelihood of reflection, we give a brief overview of the physiological processes included in the model. The voltage equation includes an inward depolarizing calcium (Ca) current, outward rectifying potassium (K) current, and a leak ($\ell$) current. The local reaction dynamics take the form
\begin{align} \label{eqn:ML_eqns}
F(V,n;\mu)=& \begin{pmatrix} f_1(V,n;\mu) \\  f_2(V,n;\mu) \end{pmatrix} =  \begin{pmatrix} -G_{Ca} m_{\infty}(V) \left(V - E_{Ca} \right) - G_K n \left(V - E_K \right) - G_{\ell} \left(V - E_{\ell} \right) \\
\epsilon \left[ \alpha(V) (1-n) - \beta(V)n \right] \end{pmatrix}
\end{align}
where $m_{\infty}(V) = 0.5 \left( 1 + \tanh \left(\frac{V - u_1}{u_2}\right)\right)$. Values of $G_{Ca}$ and $G_K$ give the maximum conductances for the two currents.

\begin{figure}[ht]
    \centering
        \includegraphics[width=0.8\textwidth]{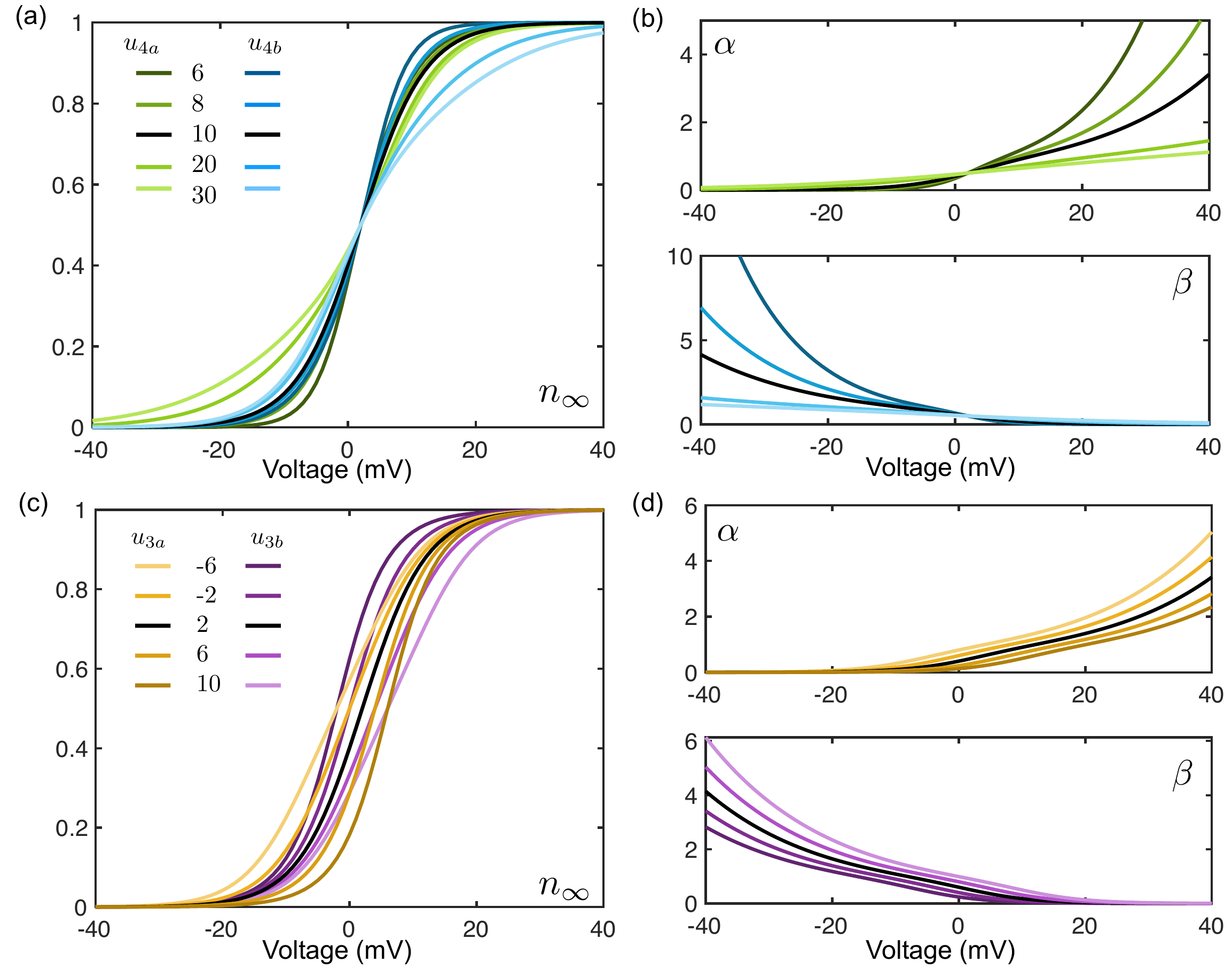}
        \caption{Impacts of $u_{3a}$, $u_{3b}$, $u_{4a}$, and $u_{4b}$ on (a,c) potassium asymptotic gating variable $n_{\infty}$, and on (b,d) activation $\alpha$ and deactivation  $\beta$ rates.  } \label{fig:ML_activation}
\end{figure}

We write the reaction equations for $n$ in terms of the opening and closing rates $\alpha(V)$ and $\beta(V)$ instead of the standard asymptotic gating variable $n_{\infty}(V)$ and time constant $\tau_n(V)$ form in order to separate the processes impacting the activation and deactivation of the potassium channels
\begin{align*}
& \alpha(V) = \frac{1}{2} \left[ \cosh \left( \frac{V - u_{3a}}{2 u_{4a}} \right) + \tanh \left( \frac{V - u_{3a}}{2 u_{4a}}\right) \cosh \left( \frac{V - u_{3a}}{2u_{4a}} \right) \right]  \nonumber\\
& \beta(V) = \frac{1}{2} \left[ 1 - \tanh \left( \frac{V - u_{3b}}{u_{4b}} \right)  \right] \cosh \left( \frac{V - u_{3b}}{2 u_{4b}} \right).\nonumber
\end{align*}
To be consistent with previous studies, the excitability parameter $\epsilon$ is factored out of $\alpha$ and $\beta$. Parameters $u_{3a}$ and $u_{4b}$ shift and scale the activation rate $\alpha$, and $u_{3b}$ and $u_{4b}$ tune the deactivation rate $\beta$. The influence of these variables on $\alpha$, $\beta$, and the asymptotic gating function $n_{\infty}$ are show in Figure~\ref{fig:ML_activation}. Increasing $u_{4a}$ or $u_{4b}$ corresponds to a reduction in the activation and deactivation rates at high and low voltages, respectively. Decreasing $u_{3a}$ shifts the activation $\alpha$ curve to lower voltages, coinciding with the opening of the potassium gates and the activation phase of the asymptotic variable $n_{\infty}$ at lower voltages. Variations to $u_{3b}$ have the same effect on the deactivation curve $\beta$, with increasing $u_{3b}$ resulting in the deactivation phase of $n_{\infty}$ shifting to higher voltages.

Throughout this study, $\epsilon$ will be the primary bifurcation parameter, and we also investigate the impacts of the parameters $\mu = \{G_{Ca}, G_K, u_{3a}, u_{3b}, u_{4a}, u_{4b}\}$. These parameters are chosen as they represent changes to the key physiological processes captured in the model, that is the depolarization and repolarization currents and the timing of the repolarization process. Parameter values are given in the Appendix, as well as the forms of $n_{\infty}(V)$ and $\tau_n(V)$.

\subsection{Traveling Pulses and Wave Trains in Reaction-Diffusion Systems}
Models of excitable media often have two traveling pulse solutions, reffered to as `fast’ and `slow’ pulses due to their speed of propagation \cite{h76,c77}. Pulses are solutions of the form
$U(x,t) = U_{f,s}(x - c_{f,s}t)$, which have a constant profile $U_{f,s}$ that is translated at speeds $c_{f,s}$ \cite{s02,kp13}; generally fast pulses are stable and slow are unstable. For rightward moving pulses, the speeds are such that $0 < c_s < c_f$ with subscripts denoting `fast’ or `slow.’ In the co-moving coordinate $\xi = x - c_{f,s}t$, the pulses are stationary solutions of 
\begin{align} \label{eqn:ss_pulse}
U_t = D U_{\xi \xi} + c_{f,s} U_{\xi} + F(U;\mu), \ \ \ \xi \in \mathbb{R}.
\end{align}

\begin{figure}[ht]
    \centering
        \includegraphics[width=0.8\textwidth]{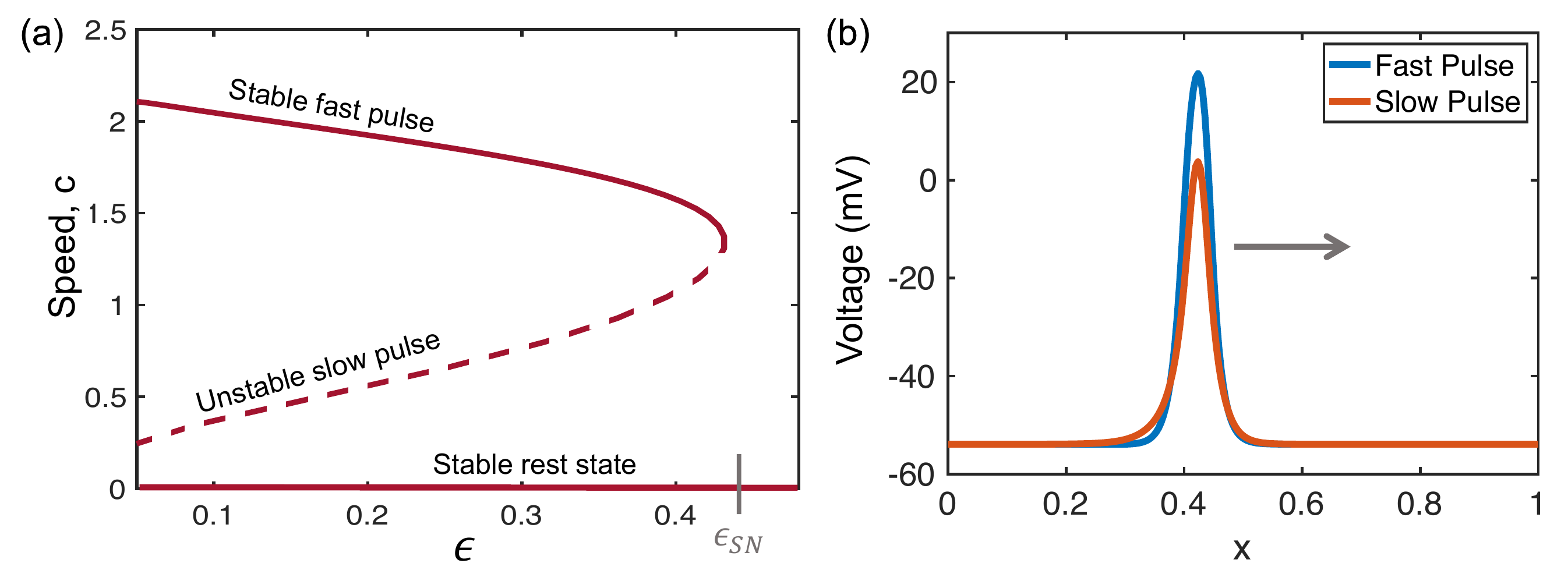}
        \caption{(a) Bifurcation diagram in $\epsilon$ for fast and slow pulses in Morris-Lecar. Default parameter set used. (b) Examples of fast and slow pulses for $\epsilon = 0.35$. Arrow indicates direction of motion. } \label{fig:fast_slow}
\end{figure}

In the Morris-Lecar system, pairs of pulses are created as the parameter $\epsilon$ moves through a saddle-node bifurcation $\epsilon_{SN}$ (Figure~\ref{fig:fast_slow}). Throughout, we will refer to $\epsilon_{SN}$ as the saddle-node bifurcation point, and we will assume the pulses exist for $0 < \epsilon < \epsilon_{SN}$. Three basic solutions exist for values of $\epsilon$ near $\epsilon_{SN}$: the fast and slow pulses, and the spatially homogeneous rest state. The stable fast pulse corresponds to action potentials that propagate in the excitable tissue and, under these conditions, the unstable slow pulse organizes the possible outcomes and acts as a threshold solution: perturbations that bring the system above the stable manifold of the slow pulse grow into a fast pulse, whereas conditions that push the system below the stable manifold will decay to the rest state.

Reaction-diffusion systems also support periodic traveling waves, i.e. wave trains, which are solutions to (\ref{eqn:rxn_diff}) of the form $U(x,t) = U_{\infty}(\kappa x - \omega t)$ where $U_{\infty}$ is $2\pi$-periodic in its argument and has wave number $\kappa$ and frequency $\omega$. Wave trains arise in one-parameter families with $\omega$ and $\kappa$ connected by the nonlinear dispersion relation $\omega = \omega(\kappa)$ \cite{s02,kp13}. In this application, periodic wave trains make up the far-field structure of the 1D spiral wave (see Figure~\ref{fig:reflection_example}b) and will be an important piece in the computation of the 1D spiral wave.

The pulse and wave train solutions are accurately computed by formulating them as equilibrium solutions in a co-moving frame \cite{amr95,dt00}. For example, the fast pulse is numerically found by solving for the vector $U_f(\xi)$ that is a solution to
\begin{align*}
0 = D U_{\xi \xi} + c_f U_{\xi} + F(U;\mu), \ \ \xi \in [0,L)
\end{align*}
on a periodic discretized domain with analytical derivatives replaced by fourth-order centered finite difference differentiation matrices. The speed $c_f$ is a free variable and is solved for with the addition of a phase condition
% \begin{align*}
% \int_0^{L} \langle U’_{old}(\xi) , U(\xi) - U_{old}(\xi) \rangle d\xi
% \end{align*}
% where the reference solution $U_{old}(\xi)$ 
that fixes the translational symmetry of the system and selects a unique solution \cite{bcdgks95}. For pulse solutions, the domain length $L$ is taken large enough so that the pulse does not interact with itself, and the domain is scaled by the wave number to be $2\pi$-periodic for wave train solutions. 

Stationary solutions can be numerically continued in any parameter to provide explicit information about how parameter changes result in modifications to the speed and pulse profile \cite{ag03, kog07}. Unlike direct numerical simulations, boundary value and continuation methods can solve for both stable and unstable patterns and saddle-node bifurcation curves. This is a crucial necessity for our study, as the slow pulses and 1D spiral waves are unstable. Throughout, we use secant methods to accurately trace out the nonlinear continuation curves \cite{bcdgks95}. 

%%%% 1D spiral
\subsection{Appearance of the one-dimensional spiral wave and reflections}

\begin{figure}[ht]
    \centering
        \includegraphics[width=0.8\textwidth]{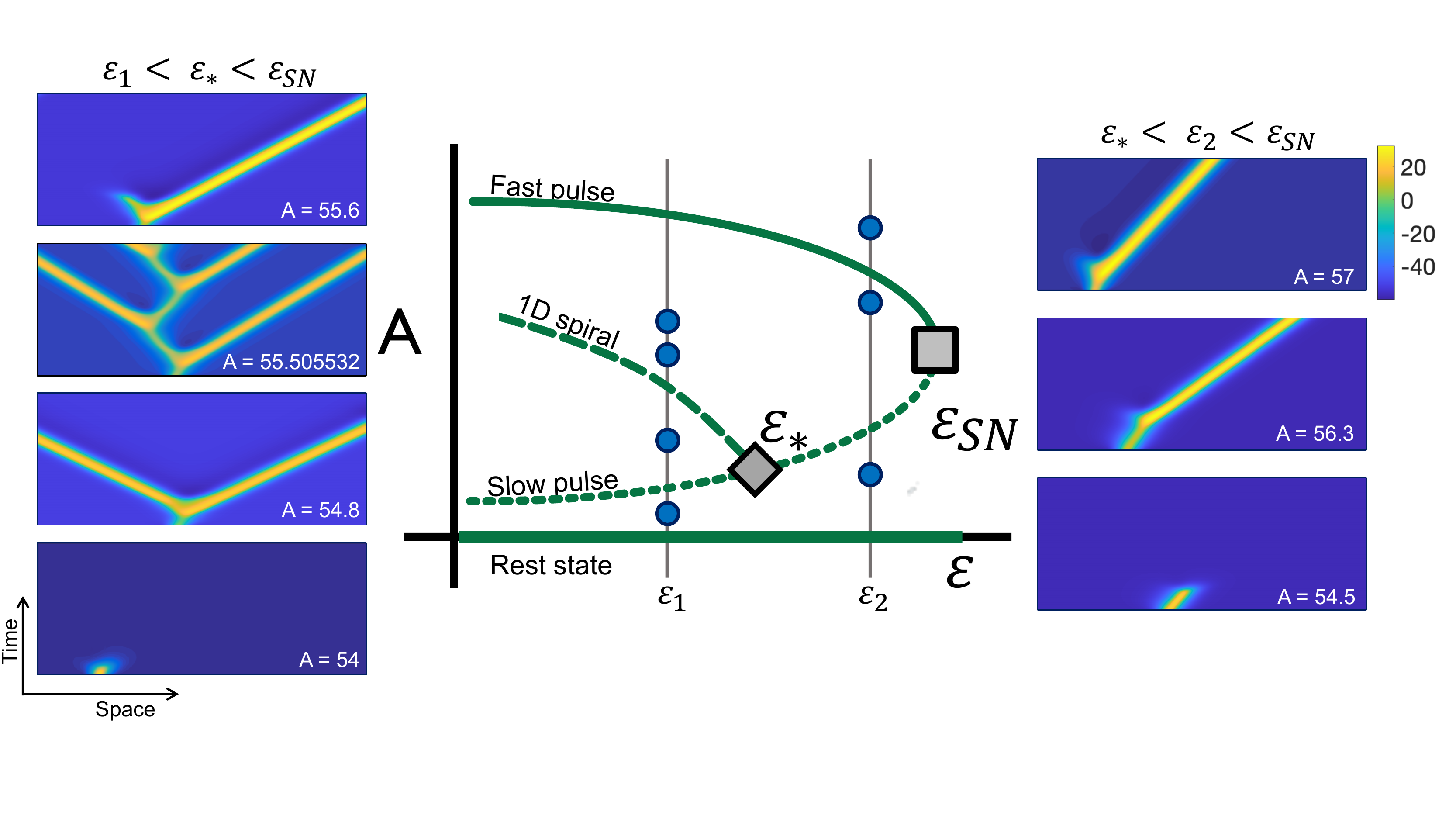}
        \caption{Schematic of bifurcation diagram based on amplitude $A$ of initial condition and parameter $\epsilon$ in the Morris-Lecar model. Locations of saddle node and heteroclinic bifurcation points indicated by $\epsilon_{SN}$ and $\epsilon_*$, respectively. Panels on left and right show resulting spatiotemporal pattern for indicated position along left and right vertical lines. Simulations run on a homogeneous domain. } \label{fig:bifur_diag}
\end{figure}

As $\epsilon$ is decreased below a critical value $\epsilon_*$, it was found \cite{cl09} that the set of solutions can become more complex as a fourth basic solution emerges: the 1D spiral wave. The 1D spiral wave is time-periodic and thus corresponds to an unstable periodic orbit \cite{er96,kh81}; throughout we will interchangeably refer to the pattern as an unstable periodic orbit or a 1D spiral wave. Even though the structure is unstable, the existence of this periodic orbit allows for new spatiotemporal dynamics, and solutions that start near or approach the orbit may lead to a transient series of reflections before diverging to a stable pattern. Figure~\ref{fig:bifur_diag} depicts the spatiotemporal patterns that can arise based on the amplitude $A$ of the initial condition and value of $\epsilon$. We emphasize that all patterns in Figure~\ref{fig:bifur_diag} are formed on a homogenous spatial domain - no heterogeneity is inducing the reflection-like behavior.

The branch of 1D spiral solutions appears to bifurcate from the branch of slow pulse solutions as illustrated in Figure~\ref{fig:bifur_diag}, and it is has been found that the period of the 1D spiral increased as $\epsilon$ approached $\epsilon_*$. This led Cytrynbaum \& Lewis to postulate that the 1D spiral arises at $\epsilon_*$ though a rearrangement of heteroclinic connections of the unstable slow pulse in a global bifurcation \cite{cl09}. 

 Our goals are to determine what conditions control the existence of the unstable periodic orbit underlying the 1D spiral pattern and factors that encourage and depress reflections when the 1D spiral does exist. The existence of 1D spiral, including the bifurcation structure giving rise it, is investigated in Section~\ref{sec:1dspiral_existence}, but we first address the likelihood of reflections in Section~\ref{sec:eigenvalues_likelihood}.
%%%%%% %%%%%% %%%%%% %%%%%% %%%%%% %%%%%% %%%%%% %%%%%% %%%%%% 
%% Reflections and the 1D spiral Wave%%%%%%%
%%%%%% %%%%%% %%%%%% %%%%%% %%%%%% %%%%%% %%%%%% %%%%%% %%%%%% 

\section{Parameter-dependent propensity for reflections} \label{sec:eigenvalues_likelihood}

First, we consider scenarios in which 1D spirals exist ($\epsilon < \epsilon_*$) and measure the parameter-dependent propensity for reflections. Specifically, we compute the 1D spiral as a spatiotemporal equilibrium solution and consider its spectral stability, as the magnitude of the unstable eigenvalue informs how repelling the periodic orbit is. The more unstable the periodic orbit, the closer a system must approach the 1D spiral's stable manifold to produce a reflection. In this section, we use direct numerical simulations to link conditions that produce reflections to the magnitude of the unstable eigenvalue. Then, we investigate parameter-dependence of the eigenvalue and determine conditions that promote or deter reflection behavior.  

%%%% Refractory experiment
\subsection{Reflections induced by elevated refractory region} \label{sec:refractory_experiment}
Reflections can be generated by a stable fast pulse interacting with a functional heterogeneity, such as a region with a transiently elevated refractory level. To demonstrate that altering system conditions influences the likelihood of reflections and to motivate the need for more precise computational methods, we test the  system's propensity for reflections using direct numerical simulation. 

\begin{figure}[ht]
    \centering
        \includegraphics[width=0.7\textwidth]{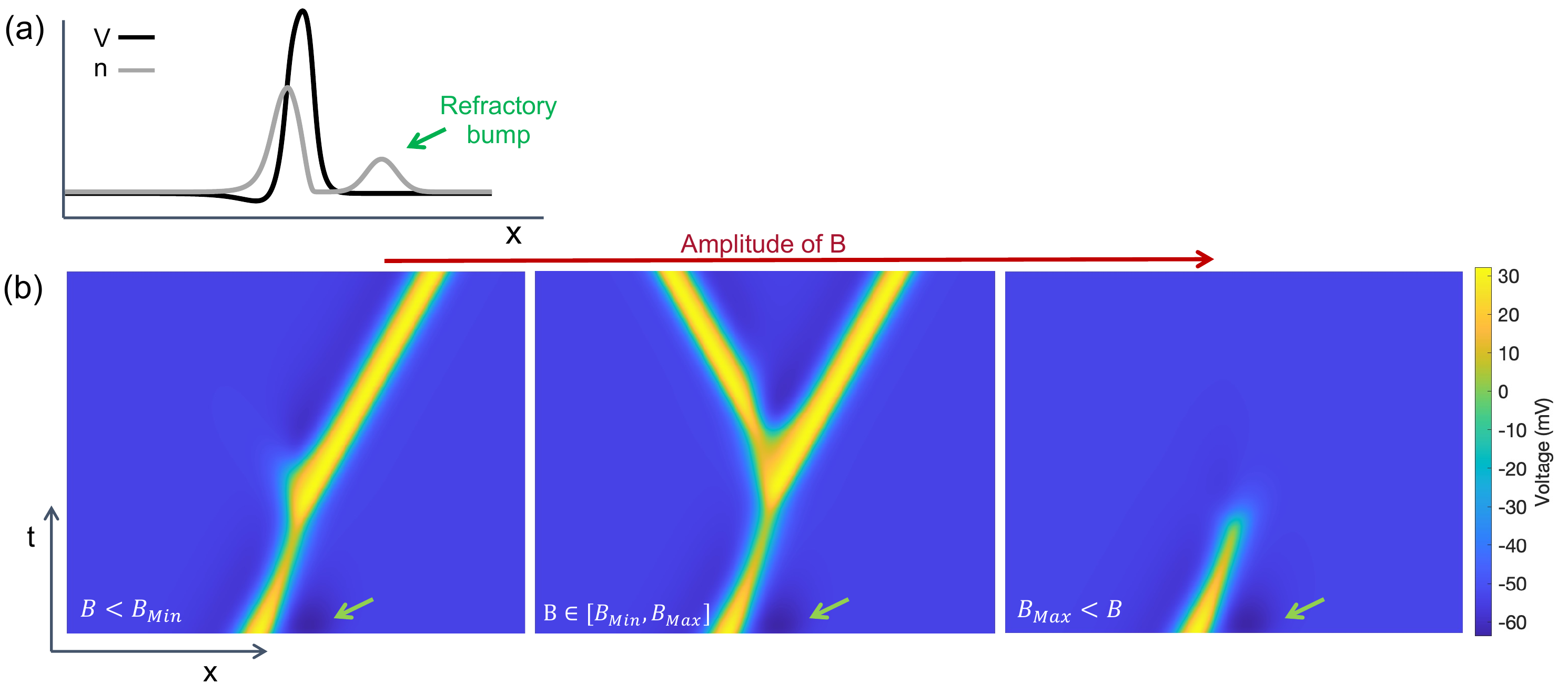}
        \caption{(a) Fast pulse interacting with a refractory bump, indicated with green arrow. (b) Series of space-time plots showing the voltage variable of a stable false pulse interacting with a refractory bump of increasing amplitude. From left to right amplitudes $B$ are 0.349, 0.34957, and 0.35.} \label{fig:reflections}
\end{figure}

The numerical experiments are performed as follows. A fast pulse is elicited so that it travels toward a small region with an elevated refractory level. First, the pulse is evolved for $t = 5$ to ensure that a stable traveling fast pulse solution is reached. A refractory bump of the form
\begin{align*}
B \exp \left( \frac{-(x - x_0)^2}{\sigma^2} \right)
\end{align*}
is added to the recovery variable $n$, where $x_0$ is set to 0.15 in front of the recovery peak and $\sigma = 0.05$. The solution is then evolved for $t = 10$, to ensure that the pulse interacts with the recovery bump. For sufficiently small amplitude $B$, refractory bumps cause the fast pulse to slow slightly, but the pulse ultimately propagates through the region and returns to the original stable form of the fast pulse. For large amplitude $B$, the pulse propagation is blocked and annihilated. Amplitudes between these two cases can generate one or more reflections (Figure~\ref{fig:reflections}). We define $ [B_{\text{min}}, B_{\text{max}}]$ as the interval of bump amplitudes that led to at least one reflection. The interval endpoints are computed via a bisection method outlined in the Appendix. In cases where reflections did occur, the length of the bump amplitude interval $B_{\text{max}} - B_{\text{min}}$ was taken to be one measure of the propensity for reflections for the given parameter set.

\begin{figure}[ht]
    \centering
        \includegraphics[width=0.5\textwidth]{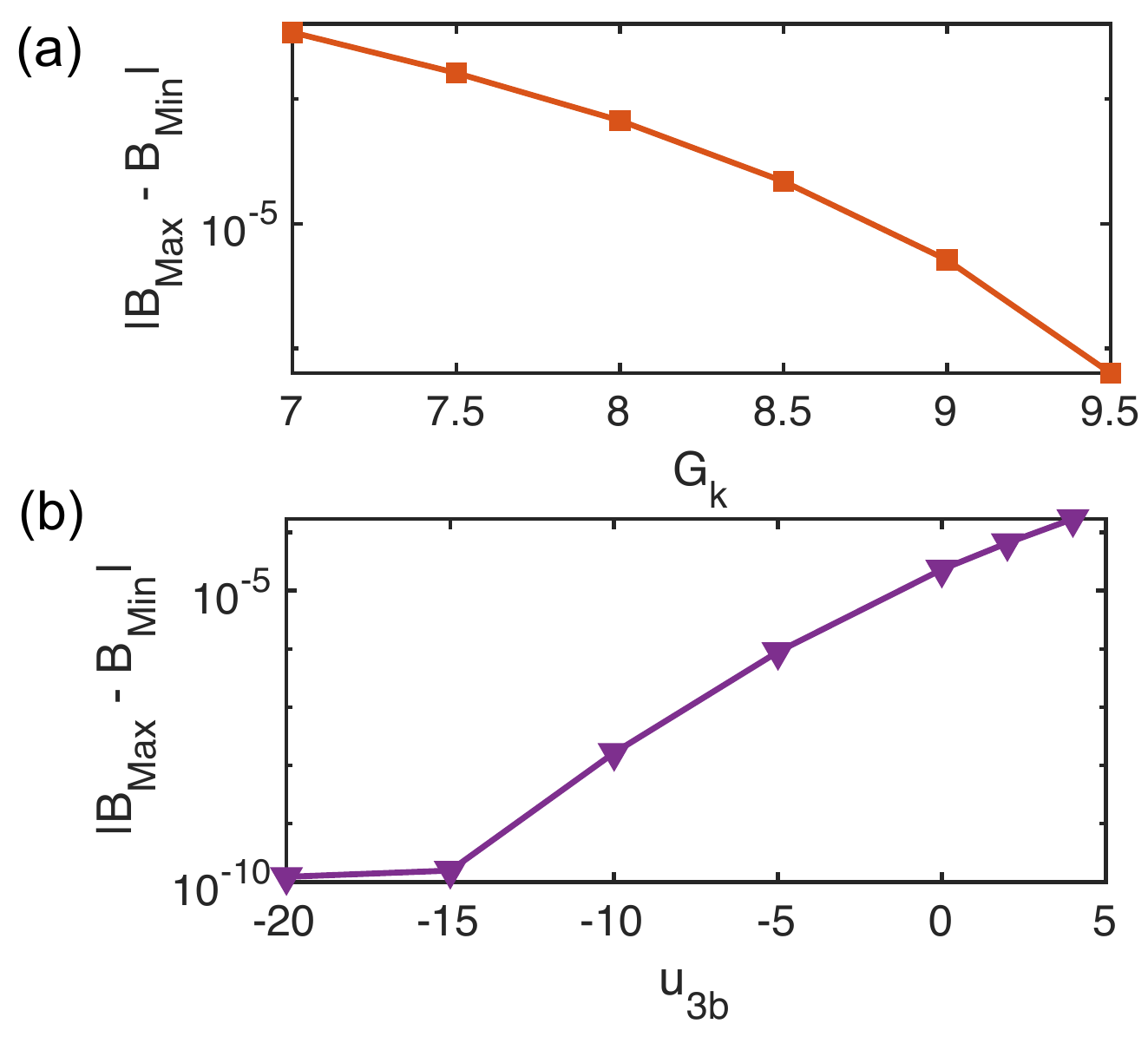}
        \caption{ Interval of refractory bump amplitudes that lead to reflection as a function of parameter (a) $G_K$ and (b) $u_{3b}$.} \label{fig:reflections_Gca_u3b}
\end{figure}

Using the experimental design described above, we test how varying the potassium conductance $G_K$ and shifting the potassium deactivation process $u_{3b}$ alters the likelihood of reflections. We find that increasing $G_K$, which corresponds to strengthening the repolarization current, results in a decreased tendency for reflections (Figure~\ref{fig:reflections_Gca_u3b}). Likewise, the bump amplitude interval and therefore propensity for reflections decreases with $u_{3b}$.

%%%% Computation of 1D spiral waves
\subsection{Computation of 1D spiral waves} \label{sec:1dspiral_comp}

The results of the refractory bump experiment yields information for how one form of heterogeneity produces reflections by perturbing the stable pulse so that the system is near the unstable 1D spiral. To efficiently determine the propensity of reflections across all types of heterogeneities and various parameters, we directly compute the 1D spiral and investigate its spectra. Specifically, we focus on the magnitude of the unstable eigenvalue because it dictates how repelling the unstable periodic orbit is; systems with a larger positive eigenvalue will diverge from the periodic orbit quicker and only initial conditions very close to the orbit will stay near it long enough to generate a transient series of reflections. Therefore, it is easier to generate reflections when the 1D spiral is less unstable (lower value of the unstable eigenvalue). To compute the 1D spiral and its stability, we consider it as a spatiotemporal pattern and explicitly take advantage of the pattern's form.  

The 1D spiral has the structure of an anti-symmetric spatiotemporal source defect \cite{ss04}: the core is a source that emits a series of counter-propagating fast pulses that assemble in the far-field into a regular structure of periodic wave trains with wave number $\kappa$ and frequency $\omega$ uniquely selected by the core. The 1D spiral waves are time $T$-periodic solutions $U_*(x,t)$ of equation~(\ref{eqn:rxn_diff}) with $U_*(x,t + T) = U_*(x,t)$. 

On the two-dimensional spatiotemporal domain $\Omega = [-L,L] \times S^1$, the 1D spiral wave can be considered as an equilibrium solution of
\begin{align} \label{eqn:1dspiral_eqn}
0 = DU_{xx} - \omega U_\tau + F(U; \mu), \ \ (x,\tau) \in [-L,L] \times S^1
\end{align}
where $\tau = \omega t$,  $\omega = 2\pi/T$, and boundary conditions at $x = \pm L$ are chosen to mimic an infinite spatial domain. By utilizing the core and far-field structure of the pattern, the 1D spiral wave can be written as
\begin{align*}
U_*(x,\tau) = W(x,\tau) + \chi(x) U_{\infty}(\kappa x - \tau)
\end{align*}
where the smooth cut-off function $\chi(x) = 1-0.5\left(\tanh[m(x + d)] - \tanh[m(x - d)]\right)$ connects the far-field wave train $U_{\infty}(\kappa x - \tau)$ with the core solution $W(x,\tau)$ \cite{ls17,gs18,ds19}. The far-field solution is pre-computed by solving for the wave train on a 1D periodic domain and interpolating onto the 2D spatiotemporal domain $\Omega$. Thus, when solving the system, the root-finding framework exploits the core and far-field structure of the pattern for added efficiency and only solves for the localized (core) function $W(x,\tau)$.

Because the 1D spiral is formulated as an equilibrium solution, its  stability can be determined from the spectra of the operator $\mathcal{L}_*$ formed by linearizing~(\ref{eqn:1dspiral_eqn}) about the 1D spiral $U_*(x,\tau)$
\begin{align*}
\mathcal{L}_*V = D V_{xx} - \omega V_{\tau} + F_U(U_*(x,\tau); \mu)
\end{align*}
and considering the resulting eigenvalue problem $\mathcal{L}_*V = \lambda V$ for $\lambda \in \mathbb{C}$. 

\begin{figure}[ht]
    \centering
        \includegraphics[width=0.5\textwidth]{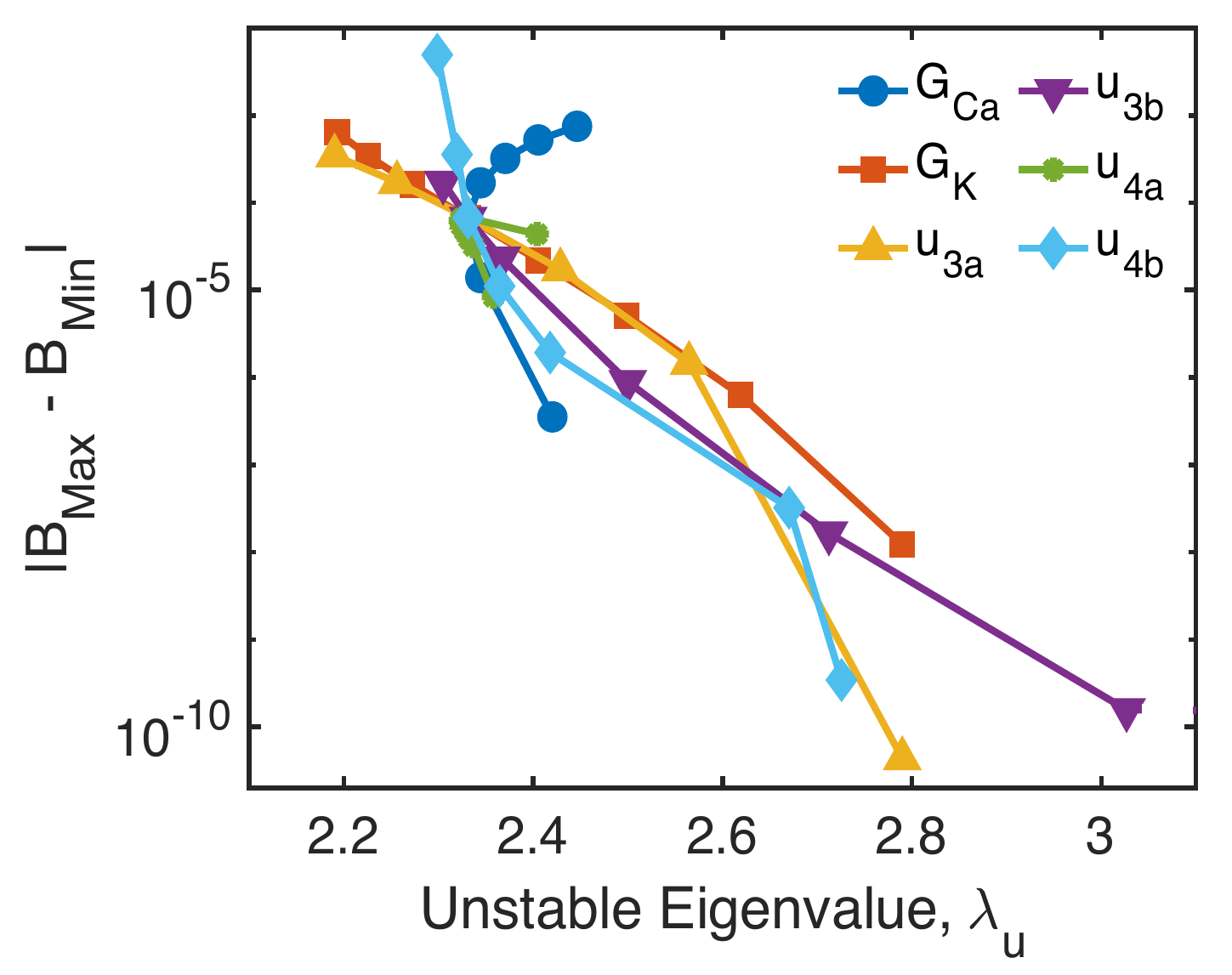}
        \caption{ Interval of refractory bump amplitudes $B$ that lead to reflections versus the unstable eigenvalue $\lambda_u$ of the 1D spiral wave. Note the log-scale of the vertical axis.} \label{fig:reflections_evals}
\end{figure}

Formulated in the above way, computing the 1D spiral as a spatiotemporal pattern corresponds to directly finding the unstable periodic orbit. Furthermore, we locate a single, real unstable eigenvalue $\lambda_u$ with corresponding eigenfunction localized near the core (shown in Appendix); the value of this unstable eigenvalue informs how repelling the orbit is. We claim reflections will be harder to generate when the eigenvalues are more positive (i.e., the system is more unstable). To support this claim, we repeat the refractory experiment for all parameters $\mu$ of interest and compare the length of the amplitude intervals that generate reflections ($B_{max} - B_{min}$) to the unstable eigenvalues of the 1D spiral. Figure~\ref{fig:reflections_evals} shows that the two quantities are correlated, corroborating the hypothesis that the 1D spiral is the fundamental dynamical mechanism underlying the ability for pulses to reflect. Thus, the unstable eigenvalue $\lambda_u$ of the 1D spiral can be used as a metric for characterizing the relative likelihood of reflections created by any heterogeneity, with more positive eigenvalues corresponding to a lower chance of reflection.

%%%%% Stability of 1D spiral and propensity for reflections
\subsection{Parameter dependent stability of 1D spiral and likelihood of reflections}

\begin{figure}[ht]
    \centering
        \includegraphics[width=0.5\textwidth]{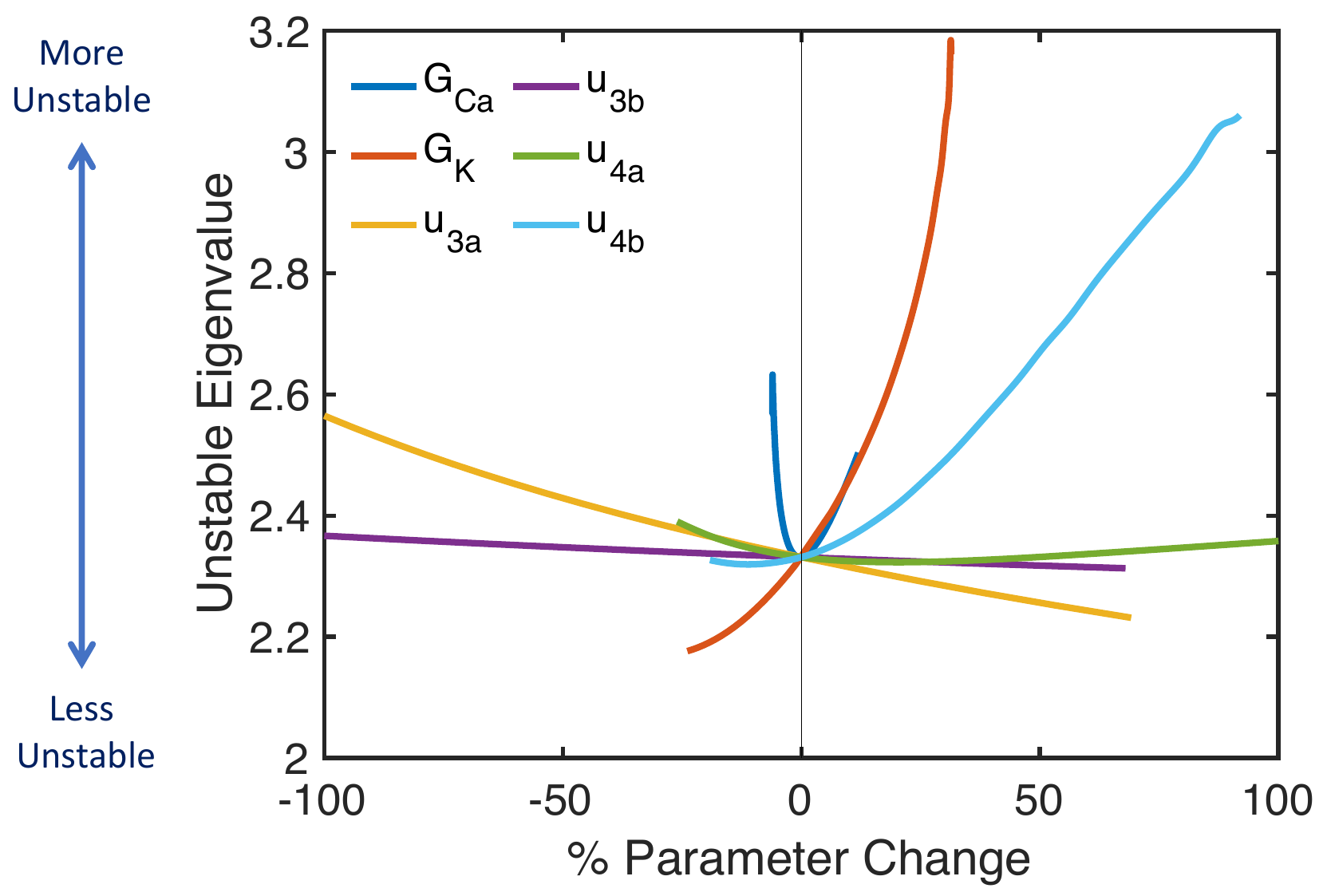}
        \caption{Results from continuation of the unstable eigenvalue $\lambda_u$ in each parameter. Displayed as a function of percent change in each parameter. Fixed $\epsilon = 0.2$. } \label{fig:unstable_eval_cont}
\end{figure}

Under the assumption that the periodic orbit exists, we seek to understand what system conditions contribute to a higher occurrence of reflection behavior. To do so, the unstable eigenvalue $\lambda_u$ is continued in parameters $\mu$ for fixed $\epsilon = 0.2$. Continuation results are reported in Figure~\ref{fig:unstable_eval_cont}. To aid in cross-parameter comparison, the eigenvalue continuation curves are shown as a function of percent change in each parameter, with the axis capped at 100\%. The unstable eigenvalue is monotonic in parameters $G_K$, $u_{3a}$, $u_{3b}$, and $u_{4b}$, and the periodic orbit becomes more unstable with decreasing $u_{3a}$ and $u_{3b}$ and increasing $G_K$ and $u_{4b}$. The eigenvalue behavior with respect to $G_{Ca}$ and $u_{4a}$ is non-monotonic, with both increases and decreases in these parameters further destabilizing the 1D spiral.

We also measure the sensitivity of the unstable eigenvalues to each parameter by (1) determining the percent change in $\lambda_u$ when $\mu$ is increased by 10\%, and (2) computing the maximum partial derivatives of the eigenvalue continuation curves with respect to the corresponding parameter over each parameter’s tested range. The percent change indicates the eigenvalue sensitivity to small changes near the default values and the maximal partial derivative provides a sensitivity index across the full tested parameter range. Both values are reported in Table~\ref{table:eval_sensitivity}; negative values in the partial derivative indicate the maximum value occurs for a decrease in the parameter value. The overall highest reduction of reflection likelihood comes from strengthening the potassium current by increasing $G_K$ ($\max \frac{\partial \lambda_u}{\partial G_{K}} = 8.06$). Combined, the two sensitivity indices indicate that the stability of the periodic orbit is most sensitive to variations in the conductances $G_{Ca}$ and $G_K$, moderately sensitive to $u_{3a}$ and $u_{4b}$, and relatively insensitive to $u_{3b}$ and $u_{4a}$. Consequences of these findings for reducing reflections and links to the physiological processes will be further discussed in Section~\ref{sec:discussion}.

\begin{table}[htp]
\caption{Parameter dependent sensitivity of the unstable eigenvalue $\lambda_u$. Column 2 gives percent change in $\lambda_u$ corresponding to a 10\% change in $\mu$. Column 3 reports the maximum value of the partial derivative of $\lambda_u$ continuation curve over each parameter’s tested range. } \label{table:eval_sensitivity}
\begin{center}
\begin{tabular}{c|c|c}
Parameter $\mu$ & \% Change $\lambda_u$ & $\max \frac{\partial \lambda_u}{\partial \mu}$\\ [5pt]
\hline 
$G_{Ca}$ & 5.69 & -2.04\\
$G_K$ & 5.46 & 8.06\\
$u_{3a}$ & -0.71 & -0.29\\
$u_{3b}$ & -0.13 & -0.07\\
$u_{4a}$ & -0.28 & -0.03\\
$u_{4b}$ & 1.42 & 0.1\\
\end{tabular}
\end{center}
\end{table}%

%%%%%%%%%%%%%%%%%%%%%%%%%%%%%%%%%%%%%%%%%%%%%%%%%%
%%% Evidence of a global bifurcation
%%%%%%%%%%%%%%%%%%%%%%%%%%%%%%%%%%%%%%%%%%%%%%%%%%
\section{Parameter-dependence for existence of the 1D spiral} \label{sec:1dspiral_existence}

Here, we investigate the bifurcation that leads to the formation of the 1D spiral wave and how the 1D spiral's existence depends on system parameters. Our focus is on uncovering mechanisms that limit the 1D spiral's existence in order to better understand conditions that hinder reflections. We start by confirming that the unstable periodic orbit does not originate via a local bifurcation of the slow pulse and follow this with a description of the proposed global bifurcation through which it arises. We then exploit the heteroclinic structure and develop a computation method to directly compute the global bifurcation point. Finally, we explore the parameter-dependence for the existence of the 1D spiral and conditions that lead to the formation of the 1D spiral are identified.  

%% Evidence of a global bifurcation
\subsection{Evidence of a global bifurcation}
The 1D spiral waves are hypothesized to arise from a global rearrangement of heteroclinic connections \cite{cl09}. To verify that the unstable periodic orbit underlying the 1D spiral does not originate through a local bifurcation from the slow pulse, we check spectrum of the slow pulse. The spectrum of a pulse consists of discrete point eigenvalues and curves of continuous spectrum; instabilities in either of these sets results in a local bifurcation. Near the saddle-node bifurcation at $\epsilon_{SN}$, the slow pulse has one real unstable point eigenvalue and a one-dimensional unstable manifold. We compute the discrete spectrum of the slow pulse to determine if any additional unstable eigenvalues appear. For $\epsilon < \epsilon_{SN}$, fast and slow pulses were computed as solutions of (\ref{eqn:ss_pulse}) and numerical continuation was used to trace out the saddle node bifurcation in $\epsilon$ (Figure~\ref{fig:fast_slow}). Eigenvalues in the discrete spectrum of the slow pulse were computed by linearizing the reaction-diffusion operator about the slow pulse \cite{amr95}. We find that spectrum of the slow pulse only contains one unstable point eigenvalue for $0 < \epsilon < \epsilon_{SN}$; the eigenvalue corresponding to the saddle-node bifurcation. 

Instabilities of the continuous spectrum are associated with instabilities of the asymptotic rest state, in this case corresponding to the resting potential. Source defects have been shown in some instances to be the result of a continuous instability \cite{ss04}, however, this is not the case in the situation presented here. The essential spectrum was computed using continuation methods described in \cite{Rademacher:2007uh} and we find that the essential spectrum of the slow pulse remains stable for all $0<\epsilon  <\epsilon_{SN}$ (eigenvalue computations in the Appendix). Because the slow pulse does not have any spectral changes, the 1D spiral cannot form as a spectral bifurcation of the slow pulse.

%%%%% Proposed global bifurcation structure
\subsection{Proposed global bifurcation structure}

\begin{figure}[ht]
    \centering
        \includegraphics[width=0.8\textwidth]{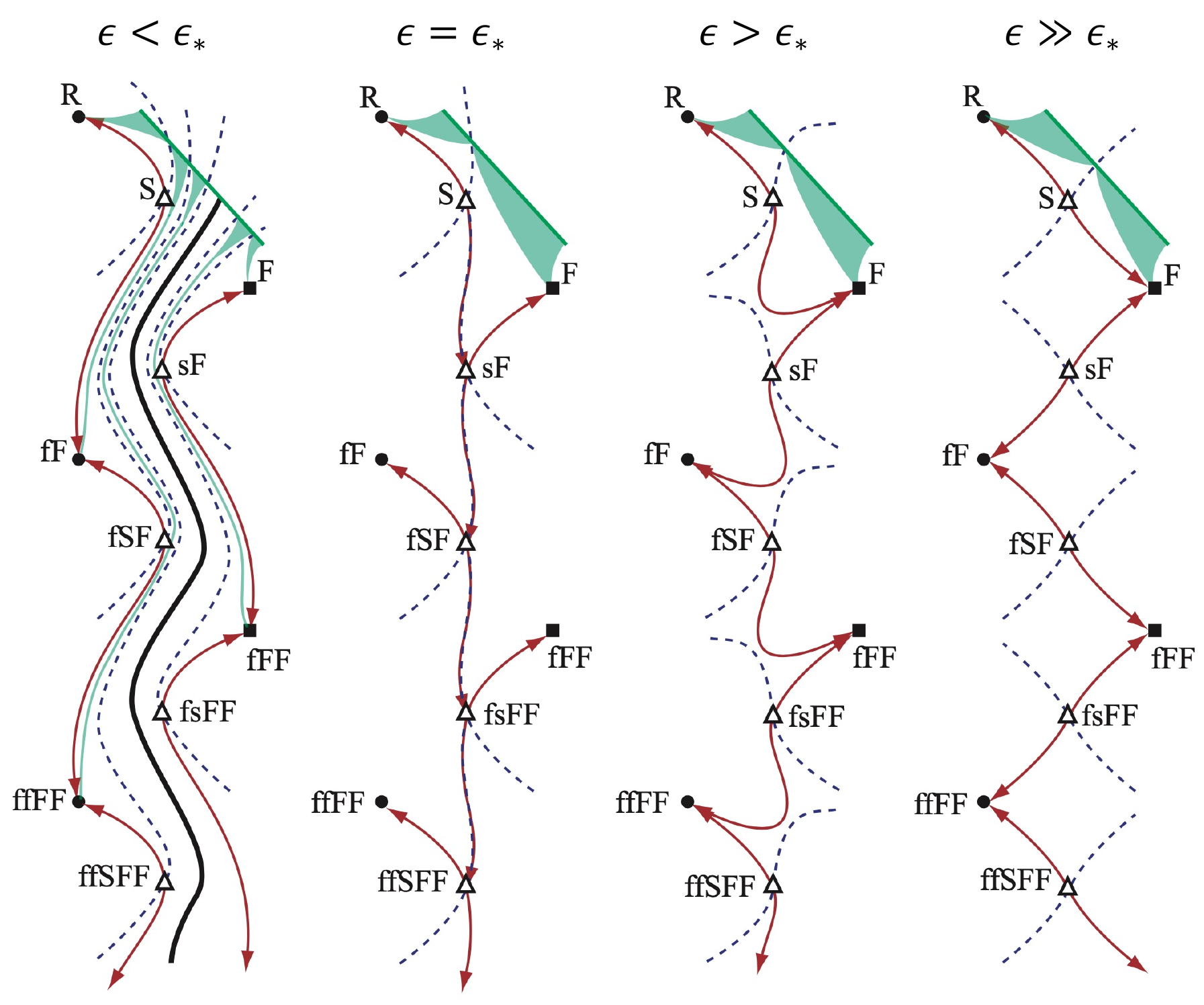}
        \caption{Schematic of proposed heteroclinic bifurcation structure as $\epsilon$ passes through $\epsilon_*$. Throughout, empty triangles represent unstable solutions; solid circles and squares are stable. For $\epsilon \gg \epsilon_*$, unstable manifolds of the slow pulse (red arrows) connect to the fast pulse $F$ and rest state $R$. At $\epsilon = \epsilon_*$, unstable manifold of the slow pulse connects to a manifold of an oppositely propagating slow-fast pulse pair (sF). Unstable periodic orbit (1D spiral) appears for $\epsilon < \epsilon_*$ (thick black curve) as unstable manifold of slow pulse connects to an oppositely propagating pair of fast pulses. Dashed lines show one-dimensional stable manifold of the slow pulse. Figure reproduced with minor revisions and permission from \cite{cl09}. } \label{fig:het_bif_struc}
\end{figure}

To introduce the mathematical setting and proposed heteroclinic bifurcation, we first summarize the behavior of a set of counter-propagating of pulses on the real line for values of $\epsilon$ near the bifurcation. Following the notation in \cite{cl09}, we denote a series of traveling pulses using $f,s$ to signify leftward fast and slow pulses, and $F,S$ to signify rightward propagating fast and slow pulses. The schematic in Figure~\ref{fig:het_bif_struc} summarizes the proposed heteroclinic bifurcation. For $\epsilon > \epsilon_*$, when perturbed, slow pulses ($S$) either decay to the rest state ($R$) or grow into a fast pulse ($F$); that is, the unstable manifold of the slow pulse $W^u(S)$ connects to $R$ and $F$. Similar behavior is seen for a slow pulse that lies in a series of pulses. For example, an initial condition near $fSF$ (left fast, right slow, right fast) will evolve into either $fF$ or $fFF$; that is, the unstable manifold of $fSF$ connects to $fF$ and $fFF$. As $\epsilon$ passes through the bifurcation point $\epsilon_*$, it is hypothesized that the unstable manifold of the slow pulse switches from being connected to the fast pulse $F$ to a pair of counter-propagating fast pulses ($fF$). Now, for initial conditions near $fSF$, the system will evolve into $fF$ if the slow pulse decays to rest, or evolve into $ffFF$ if the slow pulse splits into counter-propagating fast pulses. It is this rearrangement of unstable manifolds that gives rise to the unstable periodic orbit (thick black line in Figure~\ref{fig:het_bif_struc}). That is, the 1D spiral's formation relies on the simultaneous rearrangement of all heteroclinic connections that are linked to $W^u(S)$. 

During this rearrangement of unstable manifolds, the temporal period of the 1D spiral will limit to infinity as $\epsilon \rightarrow \epsilon_*$, and at the critical point $\epsilon_*$, the limiting behavior along $W^u(S)$ is a rightward propagating slow pulse ($S$) that splits into a rightward fast pulse and leftward slow pulse ($sF$). This behavior corresponds to a heteroclinic connection from $S$ to $sF$, or equivalently the unstable manifold of the slow pulse $W^u(S)$ to a manifold of an counter-propagating slow-fast pulse pair $M_{sF}(d)$, with solutions on the manifold parameterized by the distance $d$ between the fast-slow pulse peaks. (We will exploit this structure to compute the global bifurcation point $\epsilon_*$ in sections~\ref{sec:compbif1} and \ref{sec:compbif2}.)

\subsection{Further evidence for a heteroclinic bifurcation}
\begin{figure}[ht]
    \centering
        \includegraphics[width=0.9\textwidth]{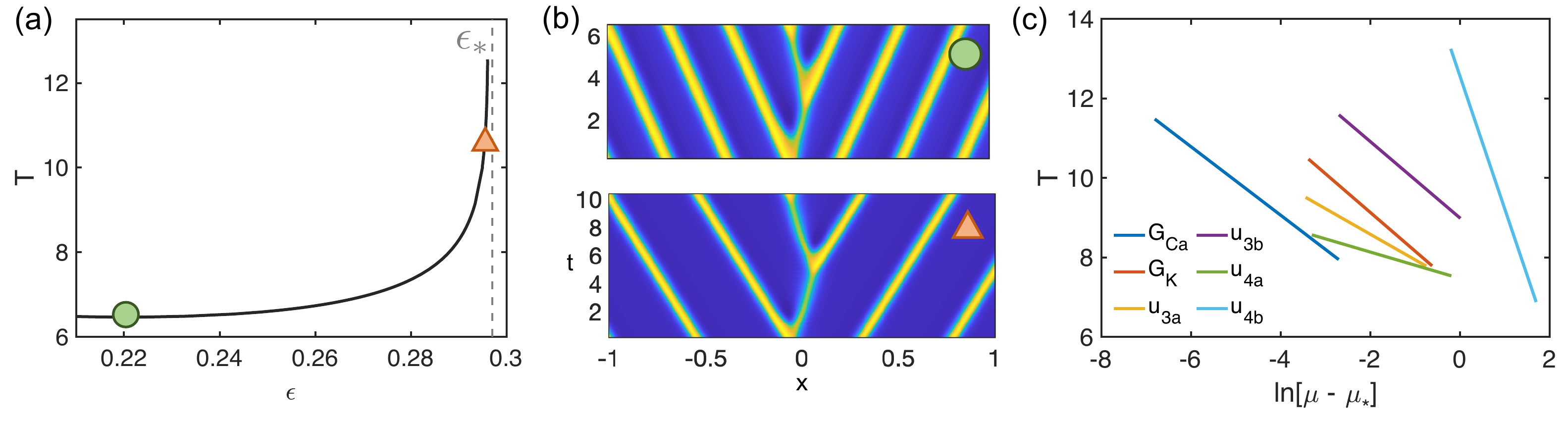}
        \caption{ Morris-Lecar Model. (a) Scaling of temporal period $T$ as solution is numerically continued in $\epsilon$ to $\epsilon_*$. (b) Spatiotemporal patterns corresponding to indicated positions on the continuation. Time $t$ rescaled to laboratory frame. (c) Temporal period $T$ versus $\log\left(\mu - \mu_*\right)$ for parameters $\mu$ near $\mu_*$.  } \label{fig:heteroclinic_bifur}
\end{figure}

To support that the global bifurcation is of the heteroclinic type, we consider the behavior of the temporal period of the 1D spiral as $\epsilon$ increases to the bifurcation point $\epsilon_*$. Using the methods described in Section~\ref{sec:1dspiral_comp}, 1D spiral waves are formulated as a spatiotemporal equilibrium pattern on the domain $\Omega = [-1,1] \times S^1$ and numerically computed. Continuation of the pattern in parameter $\epsilon$ provides information on how the temporal period $T$ and far-field spatial wave number $\kappa$ change with parameters (Figure~\ref{fig:heteroclinic_bifur}). The computations illustrate that the temporal period $T$ scales like $\log \left(\epsilon_* - \epsilon \right)$, which is the expected scaling of a heteroclinic bifurcation \cite{s00}. Consequently, as $\epsilon \rightarrow \epsilon_*$ we find that the perturbed slow pulse at the core of a 1D spiral persists for a longer time before destabilizing into a fast pulse and retrograde slow pulse, which is consistent with the proposed bifurcation structure described above.

Other system parameters $\mu$ can be used in place of $\epsilon$ when considering the 1D spiral existence interval and bifurcation point. Continuations of the 1D spiral in the parameters $\mu = \{ G_{Ca}, G_K, u_{3a}, u_{3b}, u_{4a}, u_{4b}\}$ show that the periodic orbit terminates at a point $\mu_*$, with the period $T$ scaling like $\log(\mu_* - \mu)$ for parameter values $\mu$ near $\mu_*$ (Figure~\ref{fig:heteroclinic_bifur}). These results suggest the heteroclinic bifurcation that results in the 1D spiral is a robust and more general phenomenon than just observed here.

%%%% Computation of the heteroclinic bifurcation point
\subsection{Computation of the heteroclinic bifurcation point}\label{sec:compbif1}

\begin{figure}[ht]
    \centering
        \includegraphics[width=0.85\textwidth]{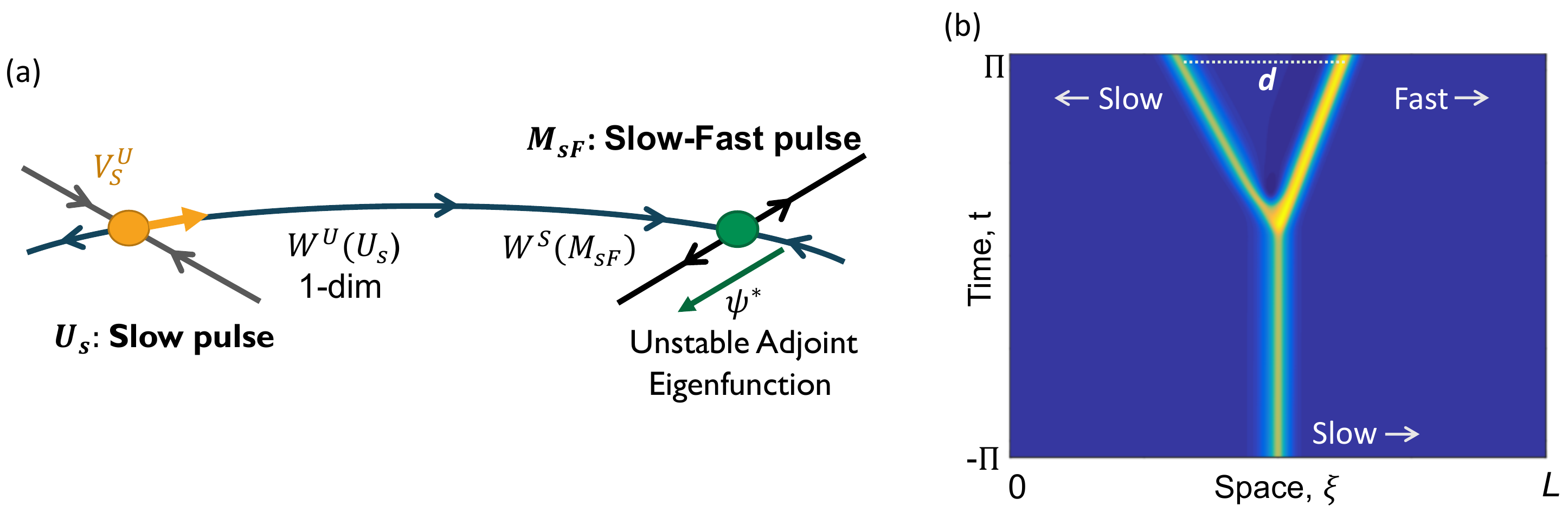}
        \caption{(a) Schematic of the hypothesized heteroclinic connection between $U_s$ and $M_{sF}$ at $\epsilon = \epsilon_*$. Truncated boundary value problem connects a solution that starts near $U_s$ to a solution that ends near $M_{sF}$. (b) Form of spatiotemporal pattern associated with the heteroclinic connection, in co-moving frame of initial slow pulse.  } \label{fig:het_bifur_alg}
\end{figure}

We develop a numerical method to directly compute the heteroclinic connection between the slow pulse ($S$) and the fast-slow counter-propagating solution ($sF$), in order to compute the global bifurcation point $\epsilon_*$ that gives rise to the 1D spiral. We provide a summary and justification for the method here, and describe implementation details in the appendix. Heteroclinic connections can be numerically computed as a truncated boundary value problem on a finite time domain by augmenting the system with projection boundary conditions that ensure the orbit approaches the equilibria in finite time \cite{bcdgks95}. We extend these methods to directly compute the heteroclinic connection between the unstable manifold of the slow pulse $W^u(U_s)$ and manifold of an oppositely propagating slow-fast pair $M_{sF}(d)$. Figure~\ref{fig:het_bifur_alg} shows a schematic of the proposed heteroclinic orbit at $\epsilon_*$ and corresponding spatiotemporal pattern that connects $S$ to $M_{sF}(d)$. Specifically, we formulate the heteroclinic connection as an equilibrium solution on a spatiotemporal domain in the co-moving frame of the slow pulse
\begin{align*}
0 = D U_{\xi\xi} + c_s U_{\xi} + F(U; \mu) - U_t, \ \ (\xi,t) \in [-L,L] \times [-T, T]
\end{align*}
equipped with the projected boundary conditions
\begin{align}
U&(\xi,t = -T) = U_s(\xi) + \gamma V_s^u(\xi) \label{eqn:cond1}\\
0 &= \langle M(\xi) - U(\xi,t = T), \psi^*(\xi) \rangle. \label{eqn:cond2}
\end{align}
Condition~(\ref{eqn:cond1}) implies that the solution at $t = -T$ is the slow pulse $U_s(\xi)$ plus a small perturbation $\gamma$ in the direction of the unstable eigenfunction $V_s^u(\xi)$. Condition~(\ref{eqn:cond2}) enforces that the difference between the solution at $t = T$ and a solution $M(\xi)$ on the manifold $M_{sF}(d)$ is perpendicular to the unstable adjoint eigenfunction $\psi^*(\xi)$ of $M(\xi)$. The necessary pulses, eigenfunctions, and adjoint eigenfunctions are additionally solved as equilibrium solutions on 1D domains, and secant continuation allows the bifurcation point $\epsilon_*$ to be continued in system parameters. 

In defining this method, we assume that the unstable manifold of $M_{sF}$ is one-dimensional, and more specifically that $W^U(M_{sF})$ corresponds to $W^u(U_s)$. This assumption is supported by results in \cite{w09} that prove systems with traveling pulse solutions also support solutions that are nearly the linear supposition of oppositely propagating pulses. Under this formulation, the eigenspaces of the counterpropagating  pulses are found to be a cross-product of the eigenspaces of each individual pulse. Because $W^U(M_{sF})$ and $W^u(U_s)$ are both one-dimensional, heteroclinic connections from $U_S$ to $M_{sF}$ only occur at isolated parameter values \cite{hs10}; for each set of parameters $\mu$, this system has one unique solution that corresponds to the global bifurcation point $\epsilon_*$. We include $\epsilon_*$ as an unknown parameter in the computations and directly solved for its value. 

\subsection{Parameter-dependent sensitivity of the Bifurcation Points \texorpdfstring{$\epsilon_*$}{} and \texorpdfstring{$\epsilon_{SN}$}{}} \label{sec:compbif2}
A major goal in this study is to assess the parameter dependent sensitivities of the heteroclinic bifurcation point $\epsilon_*$, since shifting $\epsilon_*$ will alter the existence window for reflections. However, changes to system parameters will also impact $\epsilon_{SN}$ and the existence of traveling pulses, thus we monitor the values of both $\epsilon_*$ and $\epsilon_{SN}$. Continuations of the points $\epsilon_*$ and $\epsilon_{SN}$ in $G_{Ca}$, $G_K$, $u_{3a}$, $u_{3b}$, $u_{4a}$, and $u_{4b}$ establish the independent influences of each parameter on the bifurcation points. Values of $\epsilon_*$ and $\epsilon_{SN}$ were continued in each parameter until, at minimum, either $\epsilon_*$ fell outside the range $[0.2,0.4]$ or the continuation parameter changed by at least $\pm$100\%. 

We find that $\epsilon_*$ and $\epsilon_{SN}$ change monotonically with all parameters except $u_{4a}$. Specifically, increasing $G_{Ca}$, $u_{3a}$, and $u_{3b}$ increase both bifurcation points, and increasing $G_K$ and $u_{4b}$ decreases $\epsilon_*$ and $\epsilon_{SN}$. Behaviors observed with changes in $u_{4a}$ are more complex.  Small increases in $u_{4a}$ from the default value lead to decreases in $\epsilon_*$ and $\epsilon_{SN}$, but larger increases eventually decrease both $\epsilon_*$ and $\epsilon_{SN}$; these behaviors will be further discussed below.

The sensitivity of the bifurcation points is quantified similarly to that of $\lambda_u$ in Section~\ref{sec:eigenvalues_likelihood} by computing: (1) the percent change in $\epsilon_*$ and $\epsilon_{SN}$ corresponding to 10\% increases in parameter values, and (2) the maximum partial derivative of $\epsilon_*$ and $\epsilon_{SN}$ continuations with respect to each parameter. These values present two sensitivity measures of the bifurcation points to local changes in the parameter values and are listed in Table~\ref{table:param_prcnt_change}. From these two values, we find that $\epsilon_*$ is most sensitive to changes in the calcium conductance $G_{Ca}$; this sensitivity is an order of magnitude higher than the other parameters in terms of the partial derivative. The potassium conductance $G_K$ and $u_{3a}$ have moderate influence on $\epsilon_*$, followed by $u_{4b}$, $u_{4a}$, and $u_{3b}$.

\begin{table}[htp]
\begin{center}
\caption{Parameter-dependent changes in $\epsilon_*$ and $\epsilon_{SN}$. Columns 2 and 3 are the percent change after a 10\% increase in parameter. Columns 4 and 5 are the maximum value of the partial derivative of each bifurcation point continuation over each parameter’s tested range. Negative values indicate increasing the parameter decreases the bifurcation location. }\label{table:param_prcnt_change}
\begin{tabular}{c|c|c|c|c}
Parameter $\mu$& $\%$ Change $\epsilon_*$ &  $\%$ Change $\epsilon_{SN}$ & $\max\frac{\partial \epsilon_*}{\partial_\mu}$ & $\max\frac{\partial \epsilon_{SN}}{\partial_\mu}$\\
\hline
$G_{Ca}$ & 94.36  & 63.18 & 0.827  & 0.722  \\
$G_K$    & -13.9  & -10.82 & -0.089 & -0.097   \\
$u_{3a}$ & 3.04   & 2.26 & 0.089 & 0.069    \\
$u_{3b}$ &  0.50 & 0.64 & 0.018 & 0.011   \\
$u_{4a}$ & -4.32  & 1.15 & -0.032 & 0.018   \\
$u_{4b}$ & -3.35  & -0.18  & -0.058 & -0.028  \\ 
\end{tabular}
\end{center}
\end{table}%

\begin{figure}[ht]
    \centering
        \includegraphics[width=0.9\textwidth]{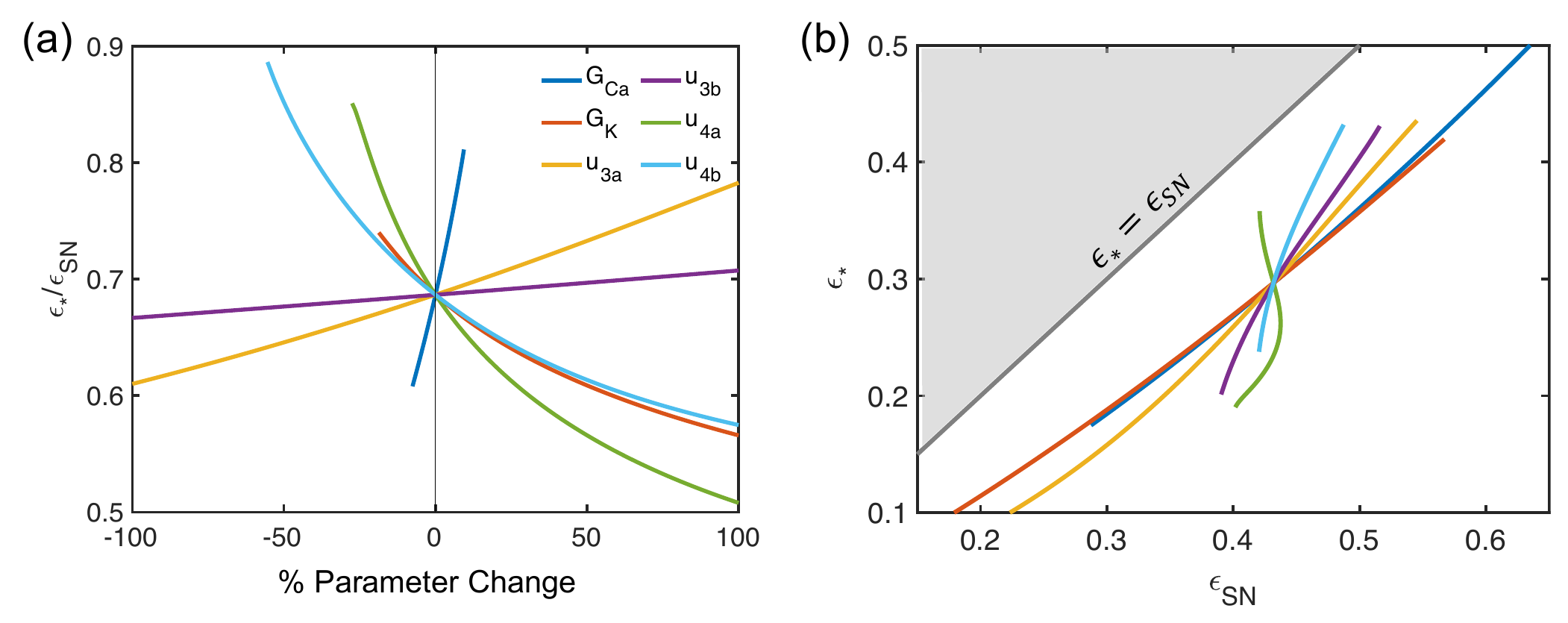}
        \caption{(a) Ratio $\epsilon_*/\epsilon_{SN}$ as function of parameter variations. (b) Relative locations of $\epsilon_*$ and $\epsilon_{SN}$.  } \label{fig:eps_sensitivity}
\end{figure}

%\subsection{Results: Relative locations of $\epsilon_*$ and $\epsilon_{SN}$} 
The ratio $\epsilon_*/\epsilon_{SN}$ provides a measure of the existence interval of the 1D spiral relative to that of the traveling pulse. Ratios near one represent conditions where reflections are possible for a large portion of the pulse existence window. The relative bifurcation values were tracked for the same changes in parameters (Figure~\ref{fig:eps_sensitivity}) and the ratio $\epsilon_*/\epsilon_{SN}$ is monotonic in all parameters. Decreases in $G_{Ca}$, $u_{3a}$, and $u_{3b}$ and increases $G_{k}$, $u_{4a}$, and $u_{4b}$ lead to reductions in $\epsilon_*/\epsilon_{SN}$. The most significant decreases correspond to increasing $u_{4a}$ and $u_{4b}$.

The ratio $\epsilon_*/\epsilon_{SN}$ does not fully characterize the relative length of the 1D spiral existence interval as both bifurcation points shift. It is also important that $\epsilon_{SN}$ remains high so that traveling pulses exist over a large range of maximal excitabilities. Corresponding values of $\epsilon_*$ and $\epsilon_{SN}$ are summarized in Figure~\ref{fig:eps_sensitivity}b. For example, while variations in $G_{Ca}$, $G_{K}$, and $u_{3a}$ reduce $\epsilon_*/\epsilon_{SN}$, these reductions in the reflection existence interval coincide with a similar reduction in the traveling pulse existence interval. These findings are consistent when comparing the sensitivities of $\epsilon_*$ and $\epsilon_{SN}$ in Table~\ref{table:param_prcnt_change}. The saddle node bifurcation point $\epsilon_{SN}$ has a similar level of sensitivity as $\epsilon_*$ to the conductances $G_{Ca}$ and $G_K$ and also $u_{3a}$. 

 More favorable outcomes (i.e., those with a decreasing ratio and relatively high $\epsilon_{SN}$) are accomplished by increases to $u_{3b}$, $u_{4a}$, and $u_{4b}$. Likewise, the sensitivity index of $\epsilon_*$ to the scaling parameters $u_{4a}$, and $u_{4b}$ is twice as great as $\epsilon_{SN}$. Moderate increases to $u_{4a}$ precipitate the most advantageous scenario: a decreasing 1D spiral existence interval coupled with a growing $\epsilon_{SN}$. Connections to the physiological processes will be revisited and discussed in Section~\ref{sec:discussion}.

%%%%%%%%%%%%%%%%%%%%%%%%%%%%%%%%%%%%%%%%%%%%%%%%%%%%%%%
%%%%%  Type I vs Type II systems
%%%%%%%%%%%%%%%%%%%%%%%%%%%%%%%%%%%%%%%%%%%%%%%%%%%%%%%
\section{Type I vs type II excitable systems} \label{sec:typeI_typeII}

In 1996, Ermentrout and Rinzel \cite{er96} used one-dimensional tissue models and systems of coupled neurons to investigate reflection behaviors. They showed that the existence of an unstable periodic orbit in a phase model of two coupled cells whose dynamics mimicked a system with a saddle-type threshold, i.e., type I excitability. From their results, they postulated that, while reflections may be possible in systems with type II excitability, they more robustly occur over a wider range of conditions in type I systems. Later results reported that reflections could not occur in type I systems under some conditions, specifically those near the traveling pulse saddle-node bifurcation \cite{cl09}, and 1D spirals have been detected in a system with type II FitzHugh-Nagumo dynamics \cite{c01}. To investigate Ermentrout and Rinzel's proposed theory for the influence of type I and type II local dynamics on the susceptibility to reflections, we interpret our results on the 1D spiral's existence and stability from Sections~\ref{sec:eigenvalues_likelihood} and \ref{sec:1dspiral_existence} in the context of the excitability type. We are particularly interested in understanding if and how type I and type II excitability influences the tendency of reflections. 

\begin{figure}[h]
    \centering
        \includegraphics[width=0.9\textwidth]{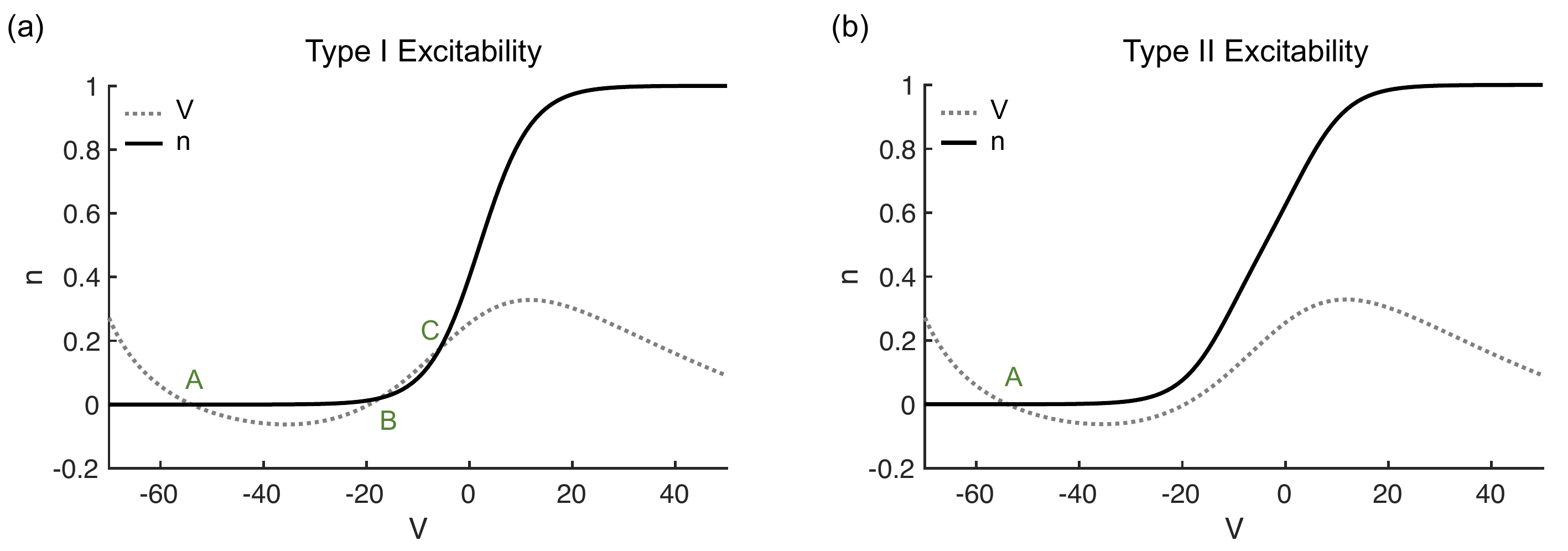}
        \caption{Nullclines in the Morris-Lecar model demonstrating (a) Type I excitability and (b) Type II excitability. Equilibrium points labeled with $A$, $B$, and $C$.  } \label{fig:excit_type_schematic}
\end{figure}

An excitable system can be classified as type I or II excitable depending on the mechanisms facilitating the transition from excitable to oscillatory dynamics in the local reaction terms and the structure of the nullclines (see Figure~\ref{fig:excit_type_schematic}). Type I systems have three equilibrium points, which we denote $A$, $B$, and $C$ in order of increasing voltage. Point $A$ is a stable node, $B$ is a saddle point, and $C$ is an unstable spiral or node. The stable manifold of $B$ acts as a strict firing threshold. Ermentrout and Rinzel postulated that it is this true threshold and subsequent consistent large amplitude action potentials characteristic of type I systems that promote robust reflection phenomena \cite{er96}. The transition from type I to type II excitability is defined by the collision and annihilation of the unstable depolarized equilibria $B$ and $C$ via a saddle-node bifurcation.  Type II systems thus have $A$ as the single equilibrium point and lack the saddle structure; these excitable systems instead display a quasi-threshold behavior. 

\begin{figure}[ht]
    \centering
        \includegraphics[width=\textwidth]{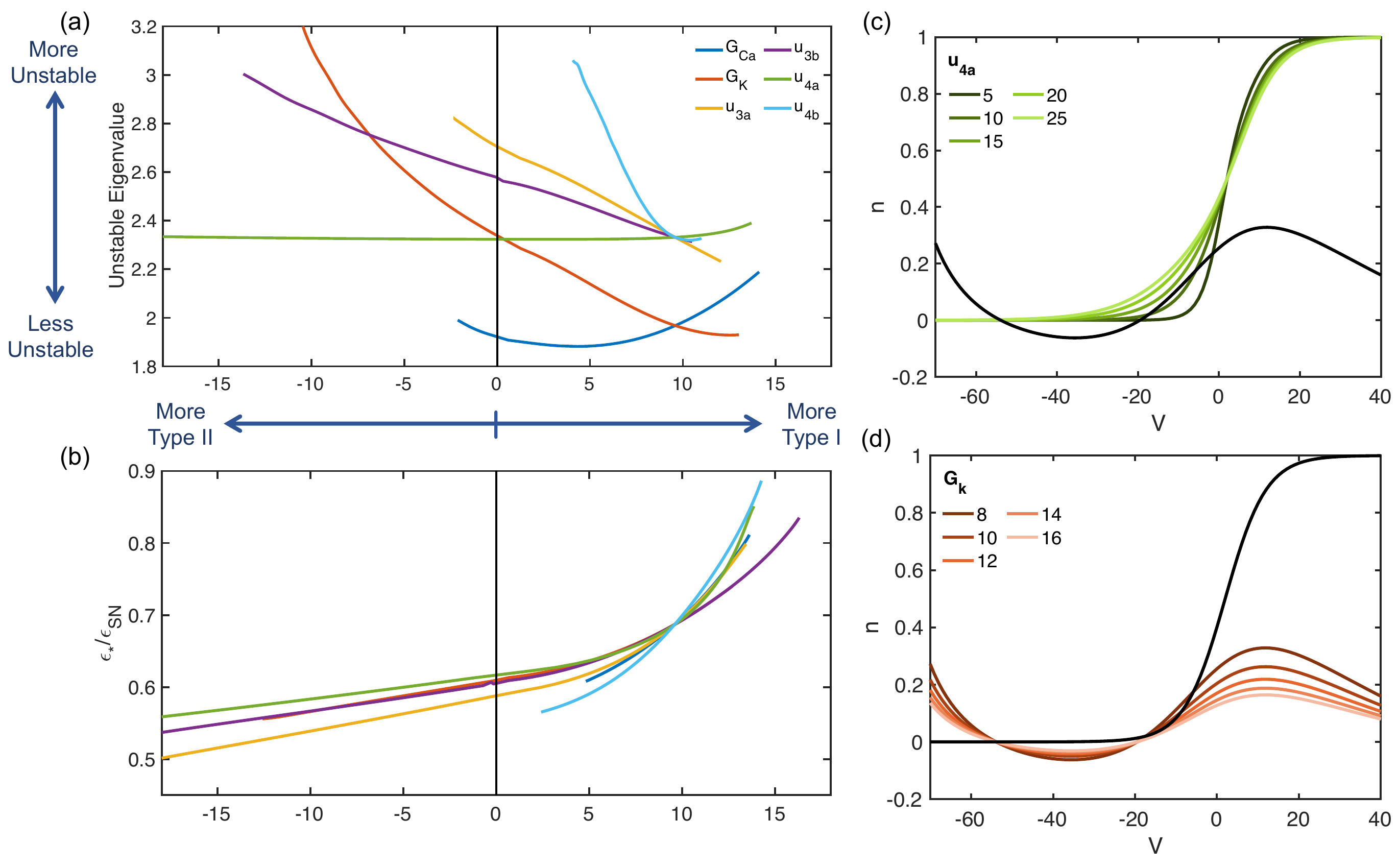}
        \caption{(a) Unstable eigenvalue continuations as a function of the type I and type II distances. Positive values on the horizontal axis indicate a larger distance between the two depolarized equilibria (b) $\epsilon_*/\epsilon_{SN}$ ratios as a function of type I and type II excitability.  (c) Phase planes with varying values of $u_{4a}$. (d) Phase planes with varying values of $G_K$.} \label{fig:exc_type}
\end{figure}

To differentiate amongst type I conditions, we classify systems with a larger Euclidean distance between the depolarized equilibria $B$ and $C$ as more ``type I-like". Likewise, type II excitable systems will be quantified by the minimum distance between the $n$- and $V$- nullclines, restricted to voltages between -30 and 0, i.e., near the $B$-$C$ saddle-node bifurcation. Conditions that lead to a larger separation in the nullcines will be considered as more ``type-II like". Figures~\ref{fig:exc_type}c-d demonstrate the transition from type I to type II excitability in a parameter that modifies the $n$-nullcline ($u_{4a}$) and one that alters the $V$-nullcline ($G_{K}$).

Figure~\ref{fig:exc_type}a plots the unstable eigenvalues as a function of the type I and type II distances. Values on the horizontal-axis represent the level of type I to type II-ness of the systems, with positive values indicating the distance between equilibria $A$ and $B$ (type I) and negative values are the distance between the $n$ and $V$ nullclines (type II). Note that with the default parameter values, the system is of type I. Eigenvalue continuations were performed with $\epsilon = 0.1$ for $G_{Ca}$ and $G_K$ to allow the 1D spiral to persist into the type II parameter range. All other eigenvalue continuations kept $\epsilon = 0.2$. Type II excitability was not achieved through significant variations in parameter $u_{4b}$ since $u_{4b}$ modifies the top bend of the $n$-nullcline with only minimal impact on the bottom bend. Likewise, due to the $\epsilon$-existence interval of the 1D spiral, eigenvalue continuations of $G_{Ca}$ and $u_{3a}$ only minimally crossed into the type II region.

Overall, the trends indicate when conditions are more type I-like, the unstable periodic orbit is less repelling, hence reflections are more likely to occur in the system. The opposite is true for type-II like systems; the further the separation in the $n$ and $V$ nullclines, the more repelling the periodic orbit is, which would result in a lower probability of reflections. Figure~\ref{fig:unstable_eval_cont} shows that parameters $G_{Ca}$ and $u_{4a}$ exhibit non-monotonic trends in the eigenvalues, yet we do find that type II local dynamics generally result in a more unstable 1D spiral. Across all parameters, the $\epsilon_*/\epsilon_{SN}$ ratio monotonically decreases as the system transitions from type I to type II excitable dynamics, again indicating reflections are less likely as they exist for a smaller portion of the traveling wave interval.

Our findings support Ermentrout and Rinzel's claim that local dynamics of type I excitability are more prone to reflections, however the results are not clear-cut. Type I systems tend to have a less repelling unstable periodic orbit, but as the eigenvalue continuations show (Figure~\ref{fig:unstable_eval_cont}) there are cases where unstable eigenvalue response is not monotonic in the shift from type I to type II. Furthermore, 1D spirals tend to exist for a larger portion of the traveling pulse existence interval in type I systems, but the decreasing $\epsilon_*/\epsilon_{SN}$ ratio often coincides with a reduced traveling wave existence interval. We furthermore find that 1D spirals exist over a broad range of parameters in type II excitable systems, and therefore reflections are possible.

%%%%%%%%%%%%%%%%%%%%%%%%%%%%%%%%%%%%%%%%%%%%%%%%%%%%%%%
%%%%%  DISCUSSION %%%%%%%%
%%%%%%%%%%%%%%%%%%%%%%%%%%%%%%%%%%%%%%%%%%%%%%%%%%%%%%%

\section{Discussion} \label{sec:discussion}

Life-threatening arrhythmia can arise when action potentials generate reflections that lead to a localized backward propagating pulse, thus it is crucial to elucidate the mathematical structure underlying reflection behavior. We utilized the Morris-Lecar model to demonstrate that distinguishing conditions that increase the propensity for generating reflections is equivalent to identifying conditions that promote the existence and stability of the unstable periodic orbit corresponding to the 1D spiral wave. 

In the Morris-Lecar system, fast and slow traveling pulses exist for values of $\epsilon$ below the saddle-node bifurcation at $\epsilon_{SN}$. For $\epsilon$ near $\epsilon_{SN}$, the slow pulse acts as the threshold separating the basins of attraction for the stable fast pulse and rest state. The possible wave behaviors become more complex for $\epsilon < \epsilon_* < \epsilon_{SN}$ as the 1D spiral wave appears. Even though the 1D spiral is an unstable structure, its existence is significant because it introduces the possibility of reflections. 

The 1D spiral wave was proposed to emerge through a  rearrangement of the unstable manifold of the slow pulse in a heteroclinic bifurcation at $\epsilon_*$ \cite{cl09}. Our findings support the proposed heteroclinic bifurcation structure. First, we confirmed the slow pulse does not undergo any spectral changes for $0<\epsilon<\epsilon_{SN}$, and additionally found that the 1D spiral period scales like $\log(|\epsilon - \epsilon_*|)$, as expected for a heteroclinic bifurcation. The central element of the heteroclinic bifurcation is that at $\epsilon = \epsilon_*$, the unstable manifold of the slow pulse rearranges and connects to a manifold of a counter-propagating slow-fast pulse pair. We directly computed this connection numerically as a truncated boundary value problem to solve for the value of $\epsilon_*$ and determined how $\epsilon_*$ changed with conditions. The computation was crucial in determining the parameter-dependence of $\epsilon_*$, but it also justifies the proposed structure. 

Numerical computations additionally provided insights into the stability of the 1D spiral. Framing the 1D spiral as a spatiotemporal equilibrium solution allowed us to compute the spectra of the linearized operator. The single positive eigenvalue of the 1D spiral informed how repelling the unstable periodic orbit was, and yielded a measure of reflection propensity.

To minimize the possibility of reflections in the Morris-Lecar model, the optimal physiological scenario is a stable fast pulse solution existing for a large range of $\epsilon$ and a highly unstable 1D spiral existing only on a narrow interval. This optimal scenario corresponds to conditions that have a high value of $\epsilon_{SN}$, a low $\epsilon_*/\epsilon_{SN}$ ratio, and a highly unstable eigenvalue for $\epsilon < \epsilon_*$. Even if optimal conditions are not likely to be fully achieved, it may be possible to suppress reflections by either shifting the system outside of the 1D spiral existence interval, or further destabilizing the periodic orbit.  

For conditions in which the 1D spiral exists, the parameters $G_{K}$, $G_{Ca}$, $u_{4b}$ and $u_{3a}$ have the most influence on the unstable eigenvalue of the 1D spiral. The magnitude of the unstable eigenvalue of the 1D spiral captures how repelling the unstable solution is, and therefore how close the system must be to the 1D spiral in order to generate a reflection. Increasing $G_{K}$ and $u_{4b}$ leads to the largest increases in the positive eigenvalue, whereas decreases in $u_{3a}$ and $G_{Ca}$ lead to more moderate eigenvalue increases. These parameter variations further destabilize the periodic orbit and reduce the likelihood of reflections as only a very narrow range of conditions will bring the system close enough to the periodic orbit to produce reflections. 

Changes in $u_{3a}$, $u_{4a}$, and $u_{4b}$ most significantly decrease $\epsilon_*$ while maintaining a large existence window of traveling pulses. Specifically, increasing $u_{4a}$ or decreasing $u_{3a}$ show large decreases in $\epsilon_*$ relative to $\epsilon_{SN}$. Modifications to the maximum calcium and potassium conductances ($G_{Ca}$ and $G_K$) greatly impact $\epsilon_*$, but also lead to similar transformations in $\epsilon_{SN}$. In fact, in these cases it is unclear if the parameter changes decrease $\epsilon_*$ directly, or if it is simply a byproduct of the decrease in $\epsilon_{SN}$. Thus, while these parameters have substantial influence on the bifurcation locations, their variations result in negative changes to the existence of the traveling pulse which, from a physiological perspective, may counteract reductions in the reflection interval. 
  
In order for pulses to propagate in excitable tissue, a balance must be struck between the strength and timing of the local excitation and recovery processes. The maximal strength of the currents is regulated by the conductances and timing via the potassium channel activation and deactivation functions. Reflections can occur when the depolarizing current is able to overcome the repolarization level at the pulse back and initiate a retrograde pulse. Therefore, it is no surprise that modifying the maximal calcium and potassium conductances considerably affects both the traveling pulse and reflection existence intervals. Furthermore, increasing the strength of repolarization ($G_K$) is the most effective at making the periodic orbit more repelling. For higher values of $G_K$, the strong repolarization current inhibits reflections as a pulse must be brought very close to the unstable periodic orbit to generate a reflection.

What is unanticipated is how altering the timing of the repolarizing potassium current disproportionately alters the ability for reflections over traveling waves. Parameters $u_{3a}$ and $u_{4b}$ favorably influence both the 1D spiral stability and existence. Decreasing $u_{3a}$ shifts the potassium activation to lower voltages, meaning the potassium channels will open  and initiate the repolarization process earlier in the action potential cycle. Increasing $u_{4b}$ scales the deactivation function $\beta$ to slower rates. Thus, the potassium channels remain open longer.

Moderate increases in $u_{4a}$ favorably shift $\epsilon_*$ to lower values while maintaining a high $\epsilon_{SN}$. These increases correspond to decreasing the activation rate $\alpha$ at higher voltages, but at threshold voltages near the peak of the slow pulse (near 0 mV), the activation rates actually increase. Therefore, the potassium channels activate more quickly earlier in the action potential and can block perturbations propagating from the back of the pulse. We suspect that further increases in $u_{4a}$ yields adverse effects and decreases $\epsilon_{SN}$, and by consequence $\epsilon_*$, as high values of $u_{4a}$ results in much slower activation rates for higher voltages. 

Based on our results in the Morris-Lecar model, we speculate that decreasing the likelihood of reflections in realistic cardiac biophysical models can be best achieved by altering the activation rates and deactivation rates of the repolarization processes. Specifically, shifting both the activation rates and deactivation rates to lower voltages may help most efficiently narrow the reflection existence interval and reduce the likelihood of producing reflections. These changes to the activation and deactivation rates act to prolong the repolarization process by opening the gates earlier in the action potential cycle and delaying their deactivation. Additionally, increasing the maximum strength of the repolarization current may also have a similar influence, but there may also be larger adverse impacts to the existence of the necessary fast pulse.

Moreover, changes that push the system to be more type-II like may also lead to a decrease in reflection behavior. Our results confirm that Ermentrout and Rinzel's hypothesis that reflections occur more robustly in systems with a saddle-threshold in the local dynamics \cite{er96}. However, we caution that while type II systems did tend to diminish reflection behavior, these systems also displayed a lower range of traveling wave existence and there was not a strictly monotonic relationship.

Furthermore, our methods described in Section~\ref{sec:1dspiral_comp} can be used to solve for any pattern that can be formulated as time-periodic source defect and investigate the underlying dynamics, stability, and existence. Reflected pulses are not unique to cardiac systems; reflections and similar patterns in the wake of traveling pulses have been observed in neurological systems \cite{Cheung:2019um, hcl76, Baccus:1998hh, bbsm00}, CO-dioxide reactions \cite{Or-GuilM2001Pbai}, and general pattern forming systems \cite{Rademacher:2004ue,ss04}. 

Our work presented here is a first step in understanding the dynamical processes that promote to the formation of reflection-induced arrhythmia and broad physiological conditions under which they occur. Moreover, several open questions remain. We introduced novel numerical methods to compute the heteroclinic bifurcation point and provided numerical evidence that the unstable periodic orbit arises through a global rearrangement of the unstable manifold of the slow pulse, but this proposition remains to be analytically proven. We speculate that the 1D spirals analyzed here and one-sided pulse generators of \cite{Yadome:2014gt} may be linked. It is likely that the two structures are connected in an unknown way in parameter space. Finally, the Morris-Lecar is an idealized model with links to realistic processes of cardiac dynamics, but further investigation of the physiological processes involved in increased reflections should be studied by applying this framework to biophysically realistic models of cardiac tissue.

%%%%%%%%%%%%%%%%%%%%%%%%%%%%%%%%%%%%%%%%%%%%%%%%%%%%%%%
%%%%% APPENDIX %%%%%%%%
%%%%%%%%%%%%%%%%%%%%%%%%%%%%%%%%%%%%%%%%%%%%%%%%%%%%%%%
\section{Appendix}

\subsection{Morris-Lecar Model}
Ranges and default values of parameters used throughout this research are included in Table~\ref{table:params_ml_full}. We write potassium channel dynamics of the Morris-Lecar model using the form of activation $\alpha$ and deactivation $\beta$ rates. These rates are related to the standard asymptotic gating variable $n_{\infty}$ and time constant $\tau_n$ by
\begin{align*}
    n_{\infty}(V) = \frac{\alpha(V)}{\alpha(V) + \beta(V)} , \ \ \tau_n(V) = \frac{1}{\alpha(V) + \beta(V)}.
\end{align*}

\begin{table}[htp]
\caption{Morris-Lecar system parameters, descriptions, default values, and tested parameter ranges. } \label{table:params_ml_full}
\begin{center}
\begin{tabular}{c|c|c| l}
\textbf{Parameter } & \textbf{Description} & \textbf{Default} & \textbf{Range}\\
\hline 
$G_{Ca}$ & Calcium conductance & 4.4 & [4.0, 4.8] \\
$G_K$ & Potassium conductance & 8 & [6, 17]\\
$u_{3a}$ & Potassium activation $\alpha$ shift & 2 & [-8.5, 4.5]\\
$u_{3b}$ &  Potassium deactivation $\beta$ shift & 2 & [-21.5, 13]\\
$u_{4a}$ &Potassium activation $\alpha$ scaling  & 10 & [7.25, 30]\\
$u_{4b}$ & Potassium deactivation $\beta$ scaling & 10 & [4.5, 22]\\
\end{tabular}
\end{center}
\end{table}%

\subsection{Numerical Methods} \label{sec:numerical_methods}
We provide implementation details of the numerical methods used and developed in this paper. All computations were performed in MATLAB 2020a \cite{MATLAB:2020}. Code to compute the 1D spiral waves, heteroclinic bifurcation, and other relevant computations can be found on GitHub \cite{code_location}.

\subsubsection{Refractory-Pulse Experiment and Bisection Method}
To create a functional heterogeneity, a Gaussian bump $g(x) = B\exp\left(\frac{-(x-x_0)^2}{\sigma^2} \right)$ is added to the recovery variable $n(x,0)$. Time evolution is carried out with the second-order IMEX Crank-Nicholason Adams-Bashforth splitting scheme on a periodic domain of length 1 discretized with 1,000 fourth-order centered finite difference spatial grid points and a temporal step size of $dt= 0.03$.

The interval of amplitudes $ [B_{\text{min}}, B_{\text{max}}]$ that lead to at least one reflection was computed with a series of bisection methods. First, an amplitude that lead to reflections $B_{\text{ref}}$ was determined by applying the bisection method to $B$ starting with values of $B = 0$ and $B = 1$. Amplitudes were classified as: less than $B_{\text{min}}$ if $V(x,t_f)$ resembled a fast pulse with one large-amplitude peak, greater than $B_{\text{max}}$ if $V(x,t_f)$ approached the homogeneous rest state, and in the interval $[B_{\text{min}}, B_{\text{max}}]$ if reflections were detected, i.e. $V(x,t_f)$ contained more than one counter-propagating large-amplitude peak. The interval bounds $B_{\text{min}}$ and $B_{\text{max}}$ were found by additional bisection methods between $B = 0$ and $B_{\text{ref}}$, and $B_{\text{ref}}$ and $B = 1$, respectively. A stopping threshold of 1e-10 was used in all bisection methods. 

\subsubsection{Computation of Traveling Pulses and Wave Trains}
Slow and fast pulses are equilibrium solutions of
\begin{align}  \label{eqn:pulse_eqn}
0 = D U_{\xi \xi} + c U_{\xi} + F(U;\mu), \ \ \xi \in [0,L)
\end{align}
in the co-moving frame $\xi = x - c t$ where $c$ is the speed of the slow or fast pulse. Pulses are solved on a periodic discretized domain with analytical derivatives replaced by fourth-order centered finite difference differentiation matrices. The speed $c_f$ is a free variable and is solved for with the addition of the phase condition
\begin{align} \label{eqn:pulse_pc}
 \int_0^{L} \langle U'_{old}(\xi) , U(\xi) - U_{old}(\xi) \rangle d\xi
\end{align}
where the reference solution $U_{old}(\xi)$ fixes the translational symmetry of the system and selects a unique solution \cite{bcdgks95}. 

Periodic wave trains are solved as stationary solutions on the co-moving $2\pi$-periodic domain $\zeta = \kappa x - \omega t$
\begin{align*}
    0 = \kappa^2 D U_{\zeta \zeta} + \omega U_{\zeta} + F(U;\mu), \ \ \zeta \in \mathcal{S}^1
\end{align*}
with wave number $\kappa$ and frequency $\omega$ connected via the nonlinear dispersion relation $\omega = \omega(\kappa)$.
A unique solution is selected by augmenting the system with the phase condition (\ref{eqn:pulse_pc}). Fourier differentiation matrices are used to approximate the derivatives. 

\subsubsection{Spectrum of Slow Pulse}
\begin{figure}[ht]
    \centering
        \includegraphics[width=0.8\textwidth]{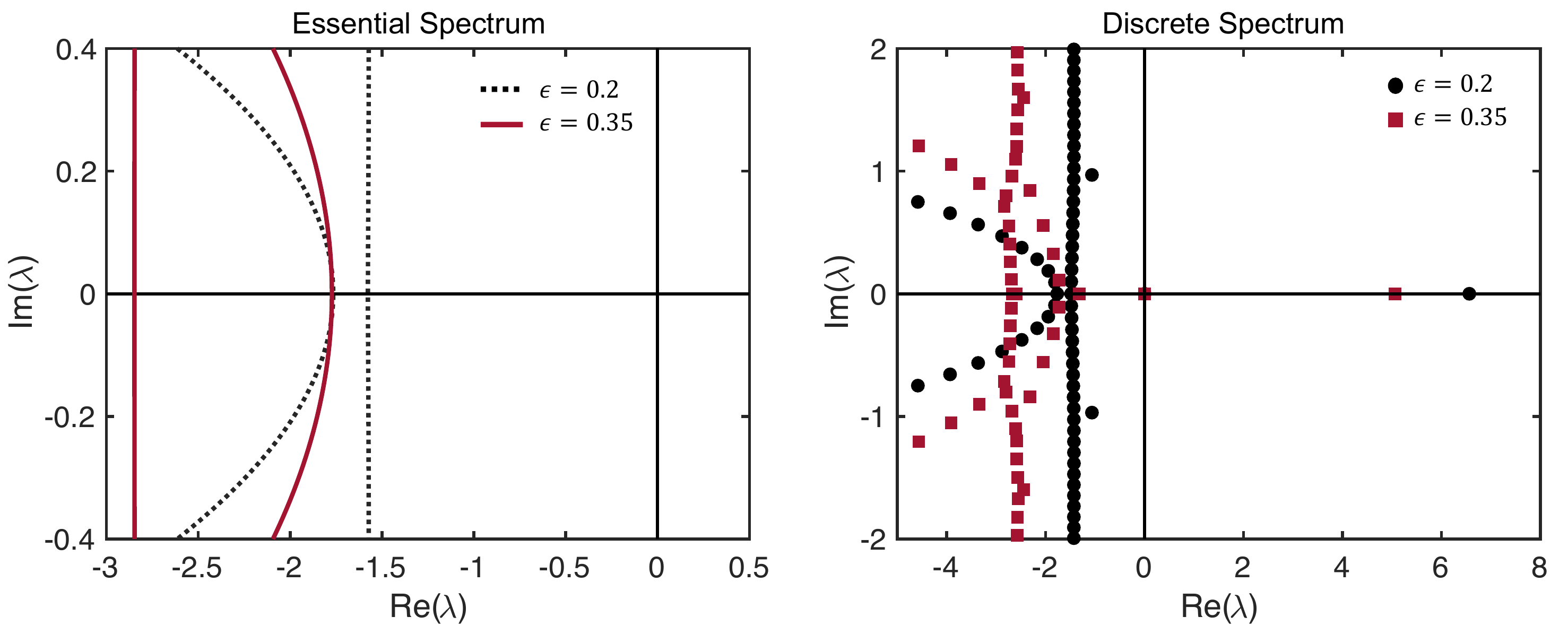}
        \caption{ Essential and discrete spectrum of the slow pulse in the co-moving frame moving at speed $c_s$.  } \label{fig:slow_pulse_spec}
\end{figure}

To confirm the 1D spiral does not emerge through a local bifurcation of the slow pulse, we check its essential and discrete spectrum. The discrete spectrum is found by linearizing the reaction-diffusion equation around the slow pulse $U_s(\xi)$ and solving for eigenvalues $\lambda$  of the  linear operator
\begin{align*}
    \mathcal{L} = D \partial_{\xi \xi} + c_s \partial_{\xi} + F_U(U_s; \mu).
\end{align*}
The essential spectrum is computed as described in \cite{Rademacher:2007uh} and provides information about the stability of the rest state at the pulse tails. As shown in Figure~\ref{fig:slow_pulse_spec}, the essential spectra is stable for values of $\epsilon = 0.2 < \epsilon_*$ and $\epsilon = 0.35 > \epsilon_*$. The discrete spectrum only has a single unstable eigenvalue for these same values, which corresponds to the one-dimensional instability generated at the saddle node bifurcation $\epsilon_{SN}$.

\subsubsection{Computation of 1D Spiral Waves}
 We consider the 1D spiral wave as a spatiotemporal equilibrium solution of 
\begin{align} \label{eqn:source_defect_app}
0 = DU_{xx} - \omega U_{\tau} + F(U;\mu), \ \ (x,\tau) \in [-L,L] \times S^1
\end{align}
with $\tau = \omega t$ and angular frequency $\omega$. Following the methods of \cite{ls17,gs18,ds19}, we split the pattern into a localized core solution $W(x,\tau)$ and the regular far-field periodic wave train pattern $U_{\infty}(\kappa x - \tau)$ by
\begin{align*}
U_*(x,\tau) = W(x,\tau) + \chi(x) U_{\infty}(\kappa x - \tau).
\end{align*}
The wave number $\kappa$ and angular frequency $\omega$ of periodic wave trains are related via the nonlinear dispersion relation $\omega = \omega(\kappa)$. Both $\omega$ and $\kappa$ are free variables to be solved for. In the implementation, the far-field periodic wave train $U_{\infty}(\kappa x - \tau)$ is first solved for on a one-dimensional $2\pi$-periodic domain with free variable $\kappa$ using Fourier spectral differentiation matrices. The periodic wave train is interpolated via Sinc interpolation onto the domain $\Omega =  [-L,L] \times S^1$ and constricted to the far-field with the cut-off function $\chi(x) = 1 - 0.5 \left(\tanh[m(x+d)] - \tanh[m(x-d)]\right)$ with $m = 10$ and $d = 0.3$.

The core function $W(x,\tau)$ is then solved for on domain $\Omega$ with free variable $\omega$ as the solution to (\ref{eqn:source_defect_app}) with Dirichlet boundary conditions $W(x = \pm L, \tau) = 0$. The phase condition
\begin{align*}
0 = \int_{2\pi/\kappa}^L \int_0^{2\pi} U'_{\infty}(\kappa x - \tau) W(x,\tau) d\tau dx
\end{align*}
matches the core and far-field solutions. 

 Fourth-order centered finite difference differentiation matrices were used in space, with Fourier spectral matrices in the temporal direction. The domain size $L = 1$ was used to allow for fast computations and allow the pattern to converge to the far-field solution. On this domain size, 401 spatial grid points and 64 temporal grid points provided more than enough accuracy for most computations. More grid points do need to be added if the temporal period increases rapidly (i.e., approaching the heteroclinic bifurcation) or for very small pulse widths. An initial condition was found via numerical time evolution and interpolated onto the spatiotemporal domain $\Omega$.

\begin{figure}[ht]
    \centering
        \includegraphics[width=0.8\textwidth]{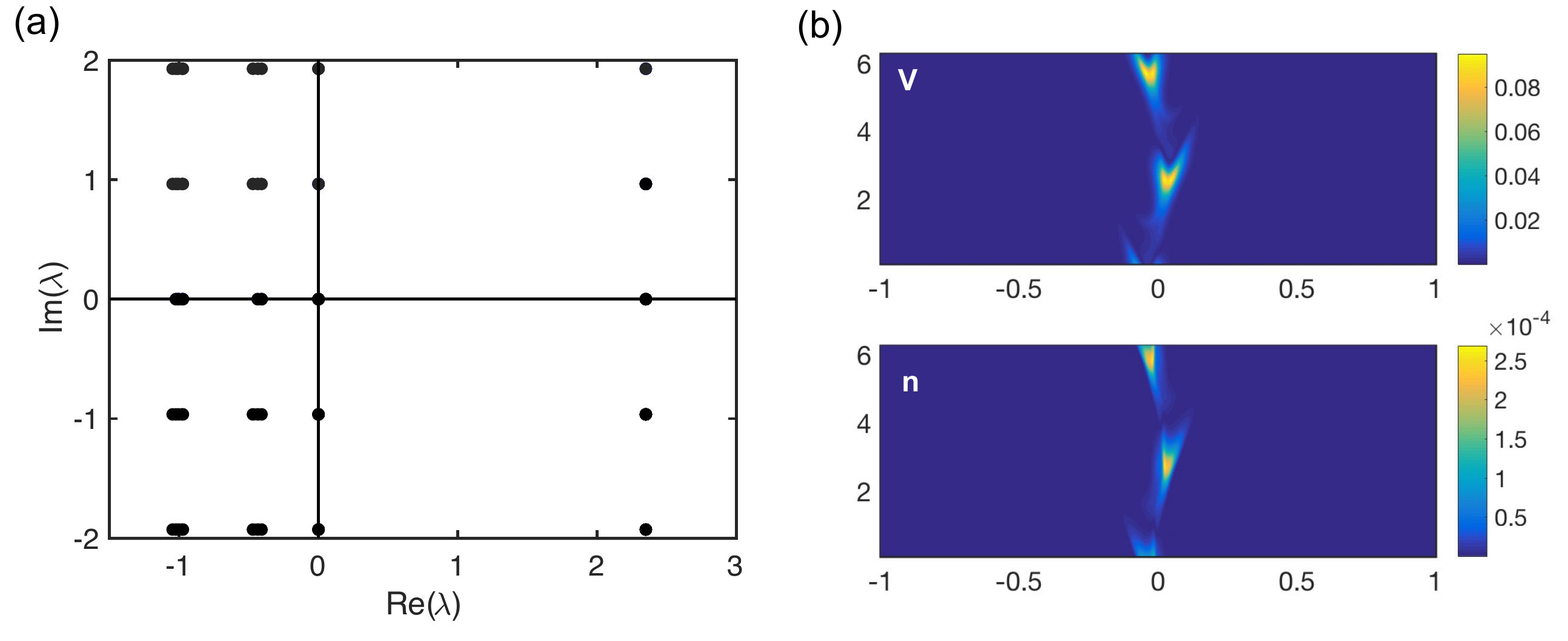}
        \caption{ (a) Eigenvalues of the 1D spiral with default parameters and $\epsilon = 0.2$. (b) Leading unstable eigenfunction, corresponding to $\lambda = 2.35$.} \label{fig:1d_sp_evals}
\end{figure}

Once a 1D spiral has been computed as a spatiotemporal equilibrium solution $U_*(x,\tau)$, the discrete eigenvalues can be numerically computed from the linearized operator
\begin{align*}
\mathcal{L} = D\partial_{xx} - \omega \partial_{\tau} + F_U(U_*(x,\tau);\mu)
\end{align*}
Leading eigenvalues were computed using the MATLAB function \texttt{eigs}. The domain size and numerical discretization are the same as those for computing the solution. The discrete eigenvalues for the default parameters with $\epsilon = 0.2$ are shown in Figure~\ref{fig:1d_sp_evals}a. Vertically periodic branches of eigenvalues arise due to the Floquet ambiguity in time. The form of the leading eigenfunction, corresponding to $\lambda = 2.35$, is shown in Figure~\ref{fig:1d_sp_evals}b. The unstable eigenfunction is localized at the core region, meaning the 1D spiral pattern destabilizes at the core. 

\subsubsection{Computation of heteroclinic bifurcation points}
 Solving for the heteroclinic bifurcation point requires computing the fast $U_f$ and slow $U_{s}$ pulses, the unstable eigenfunction of the slow pulse $V_s^u$, a oppositely propagating slow-fast pulse solution $M$ on the manifold $M_{sF}(d)$, and the unstable adjoint eigenfunction $\psi^*$ of $M_{sF}(d)$. We first describe how these components are computed, then provide the full system of equations for the heteroclinic bifurcation point. 

The fast and slow pulses are solved for on a periodic one-dimensional spatial domain as a stationary solution of (\ref{eqn:pulse_eqn}) augmented with the standard traveling wave phase condition (\ref{eqn:pulse_pc}). The unstable eigenfunction of the slow pulse $V_s^u$ solves the eigenvalue problem linearized about the slow pulse 
\begin{align}
0 &= D V_{\xi\xi} + c_s V_{\xi} + F_U(U_s;\mu)V - \lambda V \nonumber \\
0 &= \langle V, V \rangle - 1 \label{eqn:evec_pc}
\end{align}
where $\xi = x - c_s t$ is the traveling wave coordinate. Equation (\ref{eqn:evec_pc}) ensures the eigenfunction is normalized and acts as a phase condition for the eigenvalue $\lambda \in \mathbb{C}$. The 

The solution of oppositely propagating fast and slow pulses $M(\xi)$ is constructed by combining the fast and slow pulse solutions. The peak of the leftward traveling slow pulse is fixed at $\xi = q_s$, with the peak of the rightward fast pulse set to $q_f = q_s + d$. Due to the speeds of propagation of the pulses, we find that a spatial domain of $[-1,1]$ with $q_s = -0.38$ and $d = 0.64$ is large enough for the fast and slow pulses to be sufficiently separated and ensure there are no boundary effects on the computations. 

The final component is the unstable adjoint eigenfunction $\psi_*(\xi)$, which is solved for as a stationary solution of the adjoint eigenvalue problem
\begin{align*}
0 &= D \psi_{\xi \xi} - c_s \psi_{\xi} + F_U^*(U_s;\mu) \psi - \lambda \psi\\
0 &= \langle \psi, \psi \rangle -1 
\end{align*}
where the second equation normalizes the adjoint eigenfunction. In computing $\psi_*$, we assume that the unstable direction of the manifold $M_{sF}$ corresponds to the unstable direction of the slow pulse; thus the equation is that of adjoint eigenvalue problem of the slow pulse. 

To solve for the heteroclinic connection, the time is scaled by $T$ to increase efficiency and we set $\tau = t/T \in [-1,1]$. The heteroclinic connection is solved for in the co-moving frame of the slow pulse. Thus, the full system to solve for $ U(\xi,\tau) = [V(\xi,\tau), n(\xi,\tau)]^T$ on the domain $[-1,1] \times [-1,1]$ is
\begin{align}
0 &= T \left[ \delta V_{\xi\xi} + c_s V_{\xi} + f_1(V,n;\mu) \right] - V_{\tau} \label{eqn:v}\\
0 &= T \left[c_s n_{\xi} + f_2(V,n;\mu)\right] - n_{\tau} \label{eqn:n}\\
0 &= \frac{1}{T} \int_{-1}^1 \langle \partial_{\tau} V_{old}(\bar{\xi},\tau),  V(\bar{\xi},\tau) - V_{old}(\bar{\xi},\tau)  \rangle d\tau \label{eqn:temp_pc}\\
0 &= \int_0^L  \langle \partial_{\xi} V_{old}(\xi, \bar{\tau}),  V(\xi,\bar{\tau}) - V_{old}(\xi,\bar{\tau})  \rangle d\xi \label{eqn:spat_pc}\\
0 &= \langle V(\xi,\tau = 1) - M(\xi), \psi^*(\xi) \rangle. \label{eqn:t1_bc}
\end{align}
Equations~(\ref{eqn:v})-(\ref{eqn:n}) solve for the $V$ and $n$ components of the heteroclinic connection. The temporal and spatial phase conditions are given by equations (\ref{eqn:temp_pc}) - (\ref{eqn:spat_pc}).  The temporal phase condition acts at fixed $\xi = \bar{\xi}$ and spatial phase condition acts at fixed $\tau = \bar{\tau}$ with reference solution $V_{old}(\xi,\tau)$. These phase conditions are extensions of the standard traveling wave phase condition in Equation~(\ref{eqn:pulse_pc}). Equation~(\ref{eqn:t1_bc}) is the heteroclinic boundary condition at $\tau = 1$. Free parameters are $\epsilon$, $c_s$, and $T$ and are additionally solved in the system and are determined by the temporal phase condition, spatial phase condition, and heteroclinic boundary condition, respectively. 

The $\tau = -1 \ (t = -T)$ boundary condition is implemented as an external Dirichlet boundary condition. That is, the solution at $\tau = -1$ is separately defined as $U(\tau = -1,\xi) = U_s(\xi) + \gamma V^u_s(\xi)$ for a small positive $\gamma \in \mathbb{R}$, and the heteroclinic solution $U(\xi,t)$ is solved for on the remaining temporal nodes. 

The heteroclinic connection is solved for on the domain $\Omega = [-1,1] \times [-1,1]$ using fourth-order centered-finite difference differentiation matrices with 500 spatial grid points and 200 temporal grid points.

%%%%%%%%%%%%%%%%%%%%%%%%%%%%%%%%%%%%%%%%%%%%%%%%%%%%%%%%%%%%%%%%%%
\section*{Acknowledgments}
The authors express their gratitude to Bj\"{o}rn Sandstede for invaluable conversations on numerically computing the heteroclinic bifurcation point. 

\bibliographystyle{abbrv}
\bibliography{ML_bib_file}

\end{document}